\documentclass[12pt,amssymb]{amsart}

\usepackage{latexsym}
\newtheorem{propo}{Proposition}[section]
\newtheorem{defi}[propo]{Definition}

\newtheorem{lemma}[propo]{Lemma}
\newtheorem{corol}[propo]{Corollary}

\newtheorem{theo}[propo]{Theorem}

\newcommand{\ld}{,\ldots ,}

\newcommand{\ra}{ \rightarrow }

\newcommand{\lan}{ \langle }
\newcommand{\ran}{ \rangle }
\newcommand{\bl}{\begin{lemma}\label}
\newcommand{\el}{\end{lemma}}

\newcommand{\diag}{\mathop{\rm diag}\nolimits}

\newcommand{\Id}{\mathop{\rm Id}\nolimits}
\newcommand{\Irr}{\mathop{\rm Irr}\nolimits}

\newcommand{\CC}{\mathop{\bf C}\nolimits}

\newcommand{\GG}{\mathop{\bf G}\nolimits}

\newcommand{\RR}{\mathop{\bf R}\nolimits}

\newcommand{\ZZ}{\mathop{\Bbb Z}\nolimits}
\newcommand{\al}{\alpha}

\newcommand{\ep}{\varepsilon}

\newcommand{\lam}{\lambda }
\newcommand{\mar}{\marginpar }

\newcommand{\up}{^{-1}}

\newcommand{\HC}{{\mathbf H}}

\newcommand{\si}{\sigma }
\newcommand{\om}{\omega }

\newcommand{\LC}{{\mathbf L}}
\newcommand{\PC}{{\mathbf P}}
\newcommand{\GC}{{\mathbf G}}
\newcommand{\TC}{{\mathbf T}}

\def\d12{{_{12}}}

\def\acf{{algebraically closed field }}

\def\au{{automorphism }}

\def\gc{generalized character }

\def\ei{{eigenvalue }}

\def\f{{following }}

\def\ho{{homomorphism }}
\def\hos{{homomorphisms }}

\def\ii{{if and only if }}

\def\ir{{irreducible }}

\def\irr{{irreducible representation }}

\def\itf{{It follows that }}

\def\mult{{multiplicity }}

\def\po{{polynomial }}
\def\pos{{polynomials }}

\def\rep{{representation }}
\def\reps{{representations }}

\begin{document}

\title[Invariants of maximal tori 
%Quasi-projective characters for classical groups
]{
Invariants of maximal tori and unipotent constituents   of some quasi-projective characters for finite classical groups}

\author{A.E. Zalesski}
\address{Department of Physics, Mathematics and Informatics, National Academy of Sciences of Belarus, 66 Nezavisimosti prospekt, Minsk, Belarus}
\email{alexandre.zalesski@gmail.com}
\maketitle

\centerline{{\it Dedicated  to Efim Zelmanov on occasion of his 60th birthday}}

\bigskip
 \bigskip
{\small {\bf Abstract} We study the decomposition of certain reducible characters of classical groups as the sum of irreducible ones. 
Let $\GC$ be an algebraic group of classical type with defining characteristic $p>0$, $\mu$ a dominant weight and $W$ the Weyl group of $\GC$. Let $G=G(q)$ be a finite classical group, where $q$ is a $p$-power.  For a weight $\mu$ of $\GC$ the sum  $s_\mu$ of distinct  weights $w(\mu)$ with $w\in W$ viewed as a function on the semisimple elements of $G$ is known to be a generalized Brauer character of $G$ called an orbit character of $G$. We compute, for certain orbit characters and  every maximal torus $T$ of $G$, the \mult of the trivial character $1_T$ of $T$ in $s_\mu$. The main case is where $\mu=(q-1)\om$ and $\om$ is a fundamental weight of $\GC$. Let $St$ denote the Steinberg character of $G$.
Then we determine the unipotent characters occurring as constituents of   $s_\mu\cdot St$ defined to be 0   at the $p$-singular elements of $G$.
Let $\beta_\mu$ denote the Brauer character of a \rep of $SL_{n}(q)$ arising from an \irr of $\GC$ with highest weight $\mu$. 
Then we  determine the unipotent  constituents of  the characters $\beta_\mu\cdot St$ for $\mu=(q-1)\om$, and also for  some other $\mu$ (called 
strongly $q$-restricted). In addition, for strongly restricted weights $\mu$, we compute the \mult of   $1_T$ in the restriction $\beta_\mu|_T$ for 
every maximal torus $T$  of $G$.  

 }

\section{Introduction}

The groups $G$ under consideration in this papers are $GL_n(q)$, $SL_{n+1}(q)$, $Sp_{2n}(q)$, $SO_{2n+1}(q)$, $q$ odd, $SO^\pm _{2n}(q)$, $q$ odd,  $Spin^\pm_{2n}(q)$, $q$ even.  Let $\overline{F}_q$ be the algebraic closure of finite field $F_q$ of $q$ elements.
Let $\GC$ be the respective algebraic group over $\overline{F}_q$, and $W$ the Weyl group of $\GC$.  
 For the notion of a maximal torus in $G$  see \cite{DM,Ca}.
 The maximal tori of $G$, up to $G$-conjugation, are in bijection with the conjugacy classes of $W$ unless 
  $G=SO^-_{2n}(q)$, $q$ odd, and $Spin^-_{2n}(q)$, $q$ even \cite[3.3.3]{Ca}. So we denote by $T_w$ any maximal torus of $G$ from the class corresponding to $w\in W$.

Let $\TC$ be the group of diagonal matrices in $GL_{n}(\overline{F}_q)$. 
Let $\ep_i$ be the mapping sending every diagonal matrix $\diag(x_1\ld x_n)$ to the $i$-th entry $x_i$ ($1\leq i\leq n$).
There is a natural embedding $GL_{n}(\overline{F}_q)\ra \GC$ which identifies $\TC$ with a maximal torus 
of $\GC. $ So $\ep_1\ld \ep_n$  can be viewed as weights of $\GC$, as well as $\sum a_i\ep_i$ for $a_i\in \mathbb{Z}$. 
Set $\om_i=\ep_1+\cdots +\ep_i$ ($1\leq i\leq n$). 
Then $\om_i$ 
 is a fundamental weight of $\GC$, unless $i=n$ for $\GC$  of type $B_n$ and $i=n-1,n$ for type  $D_n$.
As $W$  acts on the weights of $\GC$,  we may set $W_i=\{w\in W: w(\om_i)=\om_i\}$. It is well known that
$W_i$ is the Weyl group of a certain Levi subgroup $\LC_i$ of $\GC$.  For  finite groups $A\subset B$ denote by $1_A$
the trivial character of $A$ and by $1_A^B$ the induced character. 

\begin{theo}\label{th1} Let $\GC,G,W$ be as above, $G\neq SO^-_{2n}(q),Spin^-_{2n}(q)$, and for $w\in W$ let $T_w$ be a respective maximal torus in $G$.  Let $\mu=(q-1)\om_i$, where $i\in\{1\ld n\}$. Then the number of distinct weights   $g(\mu)$ $(g\in W)$ such that  $g(\mu)(T_w)=1$
 is equal to $1_{W_i}^W(w)$. This also equals  the number of distinct weights   $g(\om_i)$ $(g\in W)$ such that $g(\om_i)(T_w)\subseteq  F_q$.
\end{theo}

Unipotent characters are introduced by Deligne and Lusztig \cite{DL}. For any character 
$\si$ of $G$ we denote by $u(\si)$ the ``unipotent part'' of $\si$, which is the sum of all unipotent \ir constituents of $\si$ regarding multiplicities.
 For the notions of  Harish-Chandra induction and the Steinberg character see \cite[Ch. 4,9]{DM}. If $\tau$ is a character of a Levi subgroup $L$ of $G$ then $\tau^{\#G}$ denotes the Harish-Chandra induced character.

\begin{theo}\label{th2} Let $\GC,G,W,\om_i, W_i$ be as above,   $\mu=(q-1)\om_i$ $(1\leq i\leq n)$ and let $s_\mu $ be the orbit character of $G$ coresponding to $\mu$.  If $G= SO^-_{2n}(q)$ or $Spin^-_{2n}(q)$, assume $j<n$. Then
$u(s_\mu\cdot St)=St_{L_i}^{\# G}$,
where $L_i$ is a Levi subgroup of $G$ with Weyl group $W_i$ and $St_i$  the Steinberg character of $L_i$.
\end{theo}

For special linear groups we have more precise results. For a dominant weight $\nu$ of 
$\GC=SL_{n+1}(\overline{F}_q)$ one can write $\nu=a_1\lam_1+\cdots +a_{n}\lam_{n}$, where $\lam_1\ld \lam_{n}$ are fundamental weights and $a_1\ld a_n$ are non-negative integers. Then $\nu$ is called $q$-restricted if 
$0\leq a_1\ld a_{n}<q$.  Let  $\nu=\nu_0+\nu_1p+\cdots +\nu_{m-1}p^{m-1}$ be the ``Steinberg expansion'' of $\nu$, where $\nu_0\ld \nu_{m-1}$  are $p$-restricted weights. We say that $\nu$ is strongly $p$-restricted if $a_1+\cdots +a_{n}< p$ and {\it strongly $q$-restricted} if each weight
$\nu_0\ld \nu_{m-1}$ is strongly $p$-restricted.

\begin{theo}\label{th5} %\mar{th5}  % and $\GC=SL_{n}(\overline{F}_q)$. 
Let   $\rho_\nu$
be an \irr of $\GC=SL_{n+1}(\overline{F}_q)$ with highest weight $\nu$ and $d_0$ the \mult of weight $0$ of $\rho_\nu$.   %Let $\nu=\nu_0+\cdots +\nu_{m-1}p^{m-1}$ be the Steinberg decomposition of $\nu$. 
Let $G=SL_{n+1}(q)$, $q=p^m$, and let $\beta_\nu$ be  the Brauer character of $\rho_\nu|_G$. 
Let $T=T_w$ be a maximal torus of $G$. Suppose that  $\nu$ is strongly $q$-restricted. 
%$\lan \nu_i,\al_0\ran\leq p-1$ for every $j=0\ld m-1$.  

$(1)$ Suppose that $\nu\neq (q-1)\lam_i$ for every $i\in \{1\ld n\}$. Then  $(\beta_\nu|_T,1_T)=d_0$ and $u(\beta_\nu\cdot St)=d_0\cdot St$. 
%, unless $\nu=(q-1)\lam_i$ for $i\in \{1\ld m-1\}$ and   $|T|\neq (q^n-1)/(q-1)$.  

$(2)$ Suppose that  $\nu=(q-1)\lam_i$ for some  $i\in\{1\ld n\}$, and let $L_i$ be a Levi subgroup of $G$
whose Weyl group is $W_i$. Then  $(\beta_\nu|_T,1_T)=d_0+1_{W_i}^W(w)$ and  $u(\beta_\nu\cdot St)=d_0\cdot St+St_{L_i}^{\# G}$, where $d_0\leq 1$. 
\end{theo}

As $\beta_\nu\cdot St$ is well known to be the  Brauer character of a projective $\overline{F}_qG$-module,
the results of this paper may be useful for  study of the decomposition matrix of the above groups
for natural characteristic $p$.  In some cases the characters $s_\mu\cdot St$ are characters of projective modules
however the question for which $\mu$ this happens in general remains open, see \cite[Ch. 9,10]{Hu}. For the cases discussed above $s_\mu\cdot St$
is expected to be a proper character vanishing at all elements of order multiple to $p$. Such characters are called quasi-projective
in \cite{WZ} and $p$-vanishing in \cite{PZ}.

\section{Deligne-Lusztig characters and $L$-controlled functions}

  {\bf Notation} Throughout this section $G$ is a finite reductive group in defining characteristic $p$,
that is, a group of
   shape $G=\GC^F$, where $\GC$ is a   connected reductive
   algebraic group, %in characteristic $p$, 
    $F$ a Frobenius morphism of $\GC$
   and $\GC^F$
stands for the set of elements fixed by $F$.   The action of $F$ on $\GC$ induces an action on every  $F$-stable
maximal  torus $\TC$, and hence on $N_{\GC}(\TC)$.
   A subgroup $T$ of $G$ is called a maximal torus if $T=\TC^F$ for
   an $F$-stable maximal torus $\TC$ of $\GC$. (In particular, a choice
   of a maximal torus $T$ means that one has also chosen $\TC$.
   The same convention is used for
    the term `a parabolic subgroup of $G$'.)  For an $F$-stable maximal torus  $\TC$ 
we set $W(T)=(N_{\GC}(\TC)/\TC)^F=N_{\GC}(\TC)^F/T$
\cite[3.13]{DM}. 

Let  $W_\GC=N_{\GC}(\TC)/\TC$  be the Weyl group of $\GC$ (we often drop the subscript and write    $W$ for  $W_\GC$ ). 
  As the group $W$ is finite,
$F$ induces an \au of  $W$. The set $\{x\in W: x\up wF(x)\}$ is called the $F$-conjugacy class
of $w\in W$, and the set $FC_{W}(w)=\{x\in W: x\up wF(w)=w\}$ is called the $F$-centralizer of $w$. 
(This is meaningful for any finite group $K$ with an \au $F$.)  The $G$-conjugacy classes of  $F$-stable maximal tori in $\GC$ are in bijection with the $F$-conjugacy classes in $W$, and we denote by $\TC_w$ a representative of  the class corresponding to  $w\in W$. We set $T_w=\TC^F_w$. 
Then  
   $W(T_w)\cong FC_{W}(w)$ \cite[3.3.6]{Ca}. The torus $\TC_1$ (that is, for $w=1$) is called the reference torus. 
The group $G$ is called non-twisted if $F$ acts
on $W$ trivially; in this case $W(T_w)\cong C_W(w)$.
To every $F$-stable maximal   torus   $\TC$
  the Deligne-Lusztig theory corresponds a generalized
character $R_{\TC,1}$.  If  $F$-stable maximal   tori $\TC $ and
$\TC'$ are $G$-conjugate then
$R_{\TC,1}=R_{\TC',1}$. An \ir character $\chi$ of $G$ is  called unipotent if $(\chi,R_{\TC,1})\neq 0$ for an $F$-stable maximal torus $\TC$. 
 
To every reductive group $\GC$ with a Frobenius morphism $F$
one corresponds the number $\ep_\GC=(-1)^r$, where $r$ is the
$F_q$-rank   of $\GC$ (see  \cite[pp.64,66]{DM}). This is
meaningful for an  $F$-stable maximal  torus $\TC$ of $\GC$.
Recall that $\ep_{\TC}=\ep_{\TC'}$ if $\TC'$ is $G$-conjugate to
$\TC$.

The notions of parabolic and Levi subgroups in $\GC$ and $G$ are standard \cite{DM}. 
(So a Levi subgroup of $G$ means a Levi subgroups of a parabolic subgroup of $G.$)
The orbit characters  $s_\nu$ are defined in the introduction, see Humphreys \cite[\S 9.6]{Hu} for more details.

Let $\phi$  be a class function on $G$. Then   $\phi |_T$ means
the restriction of $\phi$ to $T$, and, if $\theta$ is a character
of $T$ then $(\phi|_T,\theta)$ is the inner product of these
functions on  $T$. 

\subsection{$L$-controlled functions}

\begin{defi}\label{de1} %\mar{de1} 
Let $G$ be a finite reductive group and $L$ a Levi subgroup of $G.$ A function $\phi$ on G is called L-controlled if for every F-stable maximal torus $\TC$ of $\GC$ we have \begin{equation}\label{co5}
 (\phi|_T,1_T)=\sum\frac{|W(T)|}{|W_L(T')|},\end{equation}
 where the sum is over representatives of the $L$-conjugacy classes of  maximal tori $T'$  of  $L$ that are $G$-conjugate
 to $T$. (The right hand side is defined to be zero if $L$ contains no torus $G$-conjugate to $T$.) 
\end{defi}

Clearly, the values of $\phi$ at the non-semisimple elements are irrelevant for $\phi$ being $L$-controlled.
Note that an $L$-controlled function is non-zero. % as $L$ has at least one maximal torus. 

\begin{lemma}\label{ed2}%\mar{ed2} 
Let $\phi$ be a generalized character of $G$ and $St$ the Steinberg
character. %Let $\phi\cdot St=\sum_i  a_i\chi_i$, where $\chi_i\in \Irr G$ and $a_i\in \ZZ$, and  $u(\phi\cdot St)=\sum a_i\chi_i$ is a partial sum in which $\chi_i$ range over the unipotent characters of $G$.  %Suppose that $\phi$ is L-controlled. 
Then  \begin{equation}\label{eqx}u(\phi\cdot St)=\sum _{\TC  }
\frac{ (\phi|_T,1_T)}{|W(T)|} \ep_{\GC}\ep_{\TC}
R_{\TC,1}, \end{equation}

and

\begin{equation}\label{st4}(St, u(\phi\cdot St))=(St, \phi\cdot St)=  \sum _{\TC  }
\frac{ (\phi|_T,1_T)}{|W(T)|} ,\end{equation}
where in both  the sums $\TC$ ranges over representatives of all $G$-conjugacy classes of F-stable maximal tori  of $\GC$, and  $T=\TC^F$. In addition, $u(\phi\cdot St)=0$ \ii $ (\phi|_T,1_T)=0$ for every maximal torus of $G$.

\end{lemma}

Proof. Formula (3) is equivalent to  formula (2) in \cite[p.1911]{HZ2}
where one takes $s=1$. For the additional claim see \cite[Lemma 2.6]{HZ2}. 

To prove (\ref{st4}), let $\TC,\TC'$ be $F$-stable maximal tori  of $\GC$, and  $T=\TC^F,T'=\TC^F$.  The Deligne-Lusztig characters $R_{T,1},R_{T',1}$ are orthogonal if tori $\TC,\TC'$
are not $G$-conjugate, and $(R_{T,1},R_{T,1})=|W(T)|$ \cite[11.16]{DM}.  By  \cite[7.6.6]{Ca}, we have

%\medskip
%It is well known (see )  that
% The Steinberg character is  unipotent \cite[7.6.6]{Ca}; for $\phi=1_G$ we have

\begin{equation}\label{eq2} St=\sum _{\TC  }
\frac{ 1}{|W(T)|} \ep_{\GC}\ep_{\TC}
R_{\TC,1}.\end{equation}
Therefore, $St$ is a unipotent character, and hence $(St, u(\phi\cdot St))=(St, \phi\cdot St)$. So   (\ref{eq2}) and (\ref{eqx})
 yield (\ref{st4}). 

%where the sum ranges over representatives of all $G$-conjugacy classes 
%of $F$-stable maximal tori $\TC$ of $\GC$.

\begin{theo}\label{dd3}%\mar{dd3} 
 Let $\phi$ be an L-controlled function on G.  Then
$u(\phi\cdot St)=St_L^{\#G}$. \end{theo}

Proof. By Definition  \ref{de1}, $(\phi|_T,1_T)=0$ if $T$ is not $G$-conjugate to a torus in $L$, otherwise  
%for every $F$-stable maximal  torus $\TC$ in $\LC$ 
we have $\frac{(\phi|_{T},1_{T})}{|W(T)|}=\sum _{T'}\frac{1}{|W_L(T')|}$, where the sum is over the $L$-conjugacy classes 
of $F$-stable maximal   tori $\TC'\subset  \LC$ that are $G$-conjugate to $\TC$. 
Therefore, (\ref{eqx}) implies

%\medskip
% Fix $T$. Then 

\begin{equation}\label{eq22}
u(\phi\cdot St)=\sum _{\TC  }
\frac{ (\phi|_T,1_T)}{|W(T)|} \ep_{\GC}\ep_{\TC}
R_{\TC,1}= \sum _{\TC'  }
\frac{ 1}{|W(T')|} \ep_{\GC}\ep_{\TC'}
R_{\TC',1}, \end{equation}
where the right hand side sum is over representatives of the  $ L$-conjugacy classes 
of $F$-stable maximal  tori  $\TC'\subset  \LC$. 

Note that $\ep_{\GC}\ep_{\TC}=\ep_{\LC}\ep_{\TC'}$ and $R_{\TC',1}=(R^L_{\TC',1})^{\# G}$, see \cite[7.4.4]{Ca} (here $
R^L_{\TC',1}$ is the Deligne-Lusztig character of $L$).  So

$$ \sum _{\TC'  }
\frac{ 1}{|W(T')|} \ep_{\GC}\ep_{\TC'}
R_{\TC',1}=
(\sum _{\TC'  }
\frac{ 1}{|W(T')|} \ep_{\LC}\ep_{\TC'}
R^L_{\TC',1})^{\# G}.$$
By formula (\ref{eq2}), applied to $L$, $\sum _{\TC'  }
\frac{ 1}{|W(T')|} \ep_{\LC}\ep_{\TC'}
R^L_{\TC',1}=St_L$. So the result follows. 

 \bl{in8}%\mar{in8} 
 Let $ A\subset B$ be finite groups,  $F$ an \au of B such that $F( A)=A$  and $a\in A$.   Denote by $f_A^B(a)$  the number of distinct cosets $gA$ such that $aF(g)\in gA$.
 
 $(1)$  Let $D_a$ be the F-conjugacy class of a in B,
  and $a_1=a,a_2\ld a_k\in A$ be representatives of the F-conjugacy classes of elements in $D_a\cap A$.  
 Then $f_A^B(a)=\sum_{i=1}^k \frac{|FC_B(a)|}{|FC_A(a_i)|}$. 
 
 $(2)$ Let $\tilde B=B\cdot \lan F\ran$ be the semidirect product and $\tilde A=A \cdot \lan F\ran$.
 Let  $\tilde a=a\cdot F\in \tilde A$.   
Then  $f_A^B(a)=1_{\tilde A}^{\tilde B}(\tilde a)$.
 \el
 
 Proof.  (1) Recall that $FC_B(a)$ denotes the $F$-centralizer of $a$ in $B$ and similarly $FC_A(a_i)$. Set $X=\{x\in B: x\up a F(x)\in A\}$. 
 Then $|X|=|A|\cdot f_A^B(a)$.
 Furthermore, $X$ contains $FC_B(a)$, and hence $X$ is the union of cosets $xFC_B(a)$. For $c\in A\cap D_a$ let $x_c\in X$ be such that $x_c\up a F(x_c) =c$. Then all cosets $ FC_B(a)x_c$ are distinct, so $X=\cup _{c\in (A\cap D_a)}FC_B(a)x_c$. So $|X|=
 |FC_B(a)|\cdot |A\cap D_a|=|FC_B(a)|\cdot \sum _{i=1}^k\frac{|A|}{|FC_A(a_i)|} $,
 and the result follows. 
 
 (2) Let $BF=FB$ denote the coset of $B$ in $\tilde B$ containing $F$. Then $AF\subseteq BF$.
 Let $b,b'\in B$. Then $g\up bF(g)=b'$ \ii $g\up (bF)g=b'F$ (as $F(g) $ is $FgF\up$ when the \au $F$
 is realized as conjugation by $F\in \tilde B$.
 Therefore, 
 $f_A^B(a)$ equals the number of distinct cosets $gA$ $(g\in B)$ such that $g\up (aF)g\in AF.$

    Recall that  $1_{\tilde A}^{\tilde B}(\tilde a)$ is the number of distinct cosets $h\tilde A$ in $\tilde B$
    such that $ \tilde a\cdot h\tilde A=h\tilde A$. As $\tilde B=B\tilde A$,  coset representatives  can be  chosen in $B$, so this number equals the number of distinct cosets $b\tilde A$ with $b\in B$  such that $\tilde a \cdot b\tilde A=b\tilde A$, equivalently, 
    $b\up (aF) b \in \tilde A$. As every coset $BF^i$ is $\tilde B$-invariant and $aF\in BF$,
    we have $b\up (aF) b \in (BF\cap \tilde A)=AF$. In addition, the cosets $b\tilde A,b'\tilde A$ are distinct \ii so are the cosets $bA,b'A$ in $B$. So $1_{\tilde A}^{\tilde B}(\tilde a)$ equals the number of distinct cosets   $bA$ in $ B$    such that $b\up ( aF) b\in  AF $. 
 
\medskip
Let $\GC$ be a connected reductive algebraic group with Frobenius morphism $F$, $\TC_1$ a reference
$F$-stable maximal torus of $\GC$ and $W_\GC=N_{\GC}(\TC_1)/\TC_1$ the Weyl group of $\GC$. Denote by $F_1$
the \au of $W_\GC$ obtained from the action of $F$ on   $N_{\GC}(\TC_1)$. Set $\tilde W=W_\GC\cdot \lan  F_1\ran$ and
$\tilde W_{\LC}=W_{\LC} \cdot\lan F_1\ran$ for an $F$-stable Levi subgroup of $\GC$ containing $\TC_1$. In this notation we have 

 \bl{nt2}%\mar{nt2} 
 For $w\in W_\GC$ let $\TC_w\subset \GC$ be an F-stable maximal torus corresponding to $w$. Then the right hand side of  formula $(\ref{co5})$ in Definition $\ref{de1}$  coincides with 
 $f_{W_{\LC}}^{W}(w)=1_{\tilde W_{\LC}}^{\tilde W }(wF_1)$. In particular, if  $G$ is non-twisted then this coincides with $1_{W_{\LC}}^{W}(w)$. \el
 
 Proof. Let $T_w=\TC^F$ and $T'=\TC_{w'}$ for $w'\in W_\LC$. Then $|W(T)|/|W_L(T')|=|FC_W(w)|/|FC_{W_\LC}(w)|$ or 0 if $\TC$ is not $G$-conjugate to an $F$-stable maximal torus in $\LC.$   By Lemma \ref{in8}(1) with $A=W_\LC$, $B=W_\GC$, $a=w$,  the right hand side of (\ref{co5})  in Definition $\ref{de1}$ equals $f_{W_\LC}^W(w)$, and this is equal to $1_{\tilde W_{\LC}}^{\tilde W }(wF_1)$ by
  Lemma \ref{in8}(2). 
  For non-twisted groups $F_1$ is the trivial automorphism of $W$, whence the result.

\section{Maximal  tori in classical groups}

Let $\GC$ be a reductive algebraic group, $F$ a Frobenius morphism of $\GC$ and $G=\GC^F$. 
%A maximal torus in $G$ is defined to be $\TC^F$, where $\TC$ is an $F$-stable maximal torus in $\GC$. 
Let  $\TC_1$ be a maximal torus of $\GC$,   $W=N_\GC(\TC_1)/\TC_1$  the Weyl group of $\GC$ and $w\in W$. Set $D_w:=\TC_1^{w\up \circ F}=\{t\in \TC_1: w\up F(t)=t\}=\{t\in \TC_1:  F(t)=w(t)\}$. It is known that $D_w$  is conjugate to  $T_w $ in $\GC$, see \cite[the proof of Proposititon 3.3.6]{Ca}. 
 
 For our purpose we need to describe 
$D_w$ explicitly, and we do this in a way which facilitates further computations. So
we choose  representatives $w$ of the  $F$-conjugacy classes of $W$
so that the action of $w$ on $\TC_1$ and also $D_w$ to be well described. 
The choice of $w$ is called {\it canonical}, and  the  group $D_w$ %=\{t\in \TC_1:  F(t)=w(t)\}$ will be 
is called the {\it canonical $w$-torus}. 

This is done in terms  of the action of $w$ on ${\rm Hom}\,(\TC_1,\overline{F}^*_q)$, the group of algebraic group \hos from $\TC_1$ into  $\overline{F}^*_q$, the multiplicative group of $\overline{F}_q$.  The elements of ${\rm Hom}\,(\TC_1,\overline{F}^*_q)$  are called the weights of $\TC_1$ (and also of $\GC$). 
If $\dim \TC=n$ then ${\rm Hom}\,(\TC_1,\overline{F}^*_q)\cong \ZZ^n$  ($\ZZ$ is the ring of integers and $\ZZ^n$ is a free $\ZZ$-module of rank $n$). 

We first illustrate our strategy on the example of  $\GC=GL_{n}(\overline{F}_q)$.
In this case  $W\cong S_n$, where we specify $S_n$ to act 
on the set $\{1\ld n\}$. The conjugacy classes of
$W$ are parameterized by partitions $\pi=[n_1\ld n_k]$ of $n$. We fix a canonical representative  
of the conjugacy class corresponding to $\pi$  assuming  that 
$w(1)=n_1$, $w(n_1+\cdots +n_j+1)=n_1+\cdots +n_{j+1}$ for $j=2\ld k-1$ and $w(i)=i-1$ for all other numbers $i$ with $1\leq i\leq n$. 

 The reference torus $\TC_1$ of $\GC$ can be chosen to be the group of diagonal matrices. So  
every $t\in \TC_1$ can be written as   $t=\diag(x_1\ld x_n)$, $ x_1\ld x_n\in \overline{F}^\times_q$, or simply $(x_1\ld x_n)$.  The Weyl group $W$ acts on $\TC_1$ by permuting $(x_1\ld x_n)$ in a similar way, that is, $w(x_1\ld x_n)=(x_2\ld x_{n_1-1}, x_1, x_{n_1+2}\ld x_{n_1+n_2-1},x_{n_1+1}, \ldots)$ for $w\in W$.

 The mappings $\ep_i:\TC_1\ra \overline{F}_q$ given by 
 $\ep_i(t)= x_i$ ($i=1\ld n$) are sometimes called the Bourbaki weights.  Note that $\{\ep_1\ld\ep_n\}$ is a basis of $\ZZ^n$, and writing $t=(x_1\ld x_n)\in \TC_1$ is equivalent to saying that 
 $\ep_i(t)=x_i$.
 The other weights of $\TC_1$   are integral linear combinations of $\ep_1\ld \ep_n$. If $\lam=\sum z_i\ep_i$ ($z_1\ld z_n\in \ZZ$)
 then $\lam(\TC_1)=x_1^{z_1}\cdots x_n^{z_n}$.  %So $\ep_1\ld \ep_n$ form a basis in $\ZZ^n\cong {\rm Hom}\,(\TC_1,\overline{F}^\times_q)$. 
 The action of the Weyl group $W$ on   the weights of $\TC_1$ is defined by the formula $w(\lam)(t)=\lam(w(t))$ \cite[p.18]{Ca}. As $\ep_j(t)=x_j$,
 we have $w(\ep_j)(t)=\ep_j(w(t))=x_{i+1}$ unless $i=1$ or $n_1+\cdots +n_m+1$ for some $m\in\{2\ld k\}$.   Whence    $w(\ep_j)=\ep_{w\up(j)}$. In particular, $w(\ep_i)=\ep_{i+1}$ if $i\neq n_1+\cdots +n_l$ for some $l\in\{1\ld k\}$. 
 So  the action of $W$  
  on $\ep_1\ld \ep_n$ is dual to that on $\{1\ld n\}$. Observe that the $W$-orbit of $\ep_1+\cdots+\ep_j$ $(1\leq j\leq n)$  consists of weights
  $\ep_{m_1}+\cdots +\ep_{m_j}$, where $1\leq m_1<\cdots < m_j\leq n$.

\subsection{Canonical representatives of the conjugacy classes of $W$}

 In general, let  $\GC\in\{SL_{n+1}(\overline{F}_q)$,  $Sp_{2n}(\overline{F}_q), SO_{2n+1}(\overline{F}_q)$,  $q$ odd, $SO_{2n}(\overline{F}_q)$, $q$ odd,  $Spin_{2n}(\overline{F}_q)$,  $q$ even$\}$.
(If $q$  is even then $SO_{2n}(\overline{F}_q)$ is not a connected algebraic group.  So we  replace this group by  $Spin_{2n}(\overline{F}_q)$.) Furthermore, there is an injective algebraic group \ho $GL_{n}(\overline{F}_q)\ra \GC$, which identifies $\TC_1$ with a maximal torus in $\GC$.  This defines the weights $\ep_1\ld \ep_n$ for the above groups $\GC$.

Suppose first  that $\GC=SO_{2n+1}(\overline{F}_q)$ with $q$ odd, or $Sp_{2n}(\overline{F}_q)$. The Weyl group $W$ of $\GC$ in both the cases is isomorphic to the semidirect product of $S_n$ with abelian normal subgroup of order $2^n$ acting on $\TC_1$ by sending $(x_1\ld x_n)$
to $(x_1^{\pm 1}\ld x_n^{\pm 1})$. Therefore, $W$ acts on $\ZZ^n\cong {\rm Hom}\,(\TC_1,\overline{F}^\times_q)$ by sending every $\ep_i$ to $\ep_j$ or $-\ep_j$ for some $j$.

It is well known \cite{Ca1} that
 the conjugacy classes of $W$ are parameterized by
 double partitions  $\pi=[n_1\ld n_k,n_{k+1}\ld n_{k+l}]$, where $n_1\geq n_2\geq \cdots \geq n_k$, 
 $n_{k+1}\geq n_{k+2}\geq \cdots \geq n_{k+l}$ and $n_1+\cdots +n_{k+l}=n$. To avoid confusion we shall write 
$\pi=[n_1\ld n_k,n^*_{k+1}\ld n^*_{k+l}]$. We allow $k=0$ or $l=0$, in these cases we write $[- ,n^*_{k+1}\ld n^*_{k+l}]$ and 
$[n_{1}\ld n_{k}, -]$, respectively.   The canonical form of an element $w\in W$ labelled by $\pi$ 
can be described as follows. 

Let $t=(x_1\ld x_n)\in \TC_1$ and set $n':=n_1+\cdots +n_k$. Then $w$ simply permutes $x_1\ld x_{n'}$ in the way described above for $GL_{n}(\overline{F}_q)$. (For instance $x_i$ ($i\leq n'$) goes to the $(i-1)$-th position for $1<i\leq n_1$ etc.) Let $n'<i$. Then  $w$ puts $t_i$ on the $(i-1)$-th position provided $i\neq  n'+n^*_{k+1}+\cdots +n^*_{k+j}+1$ %$<i\leq n'+n^*_{k+1}+\cdots +n^*_{k+j+1}$ 
for some $j<l$. Otherwise, if $i= n'+n^*_{k+1}+\cdots +n^*_{k+j}+1$ 
for some $j<l$, then  $w$ puts $t_i\up$ on the $ (  n'+n^*_{k+1}+\cdots +n^*_{k+j+1})$-th position.  (For instance, if the double partition is $[2,3^*]$ then $w$ sends $(x_1,x_2,x_3,x_4,x_5)$ to $(x_2,x_1,x_4,x_5,x_3\up)$.)  
 
 \itf the dual action of $W$ on Hom$(\TC_1,\overline{F}_q)$ preserves the set $\{\pm \ep_1\ld \pm \ep_n\}$.   In particular,
 $w$  permutes $\ep_1\ld \ep_{n'}$ by the rule described for the  $GL_{n}(\overline{F}_q)$-case. Let $i>n'$ and let $r(j)=n_1+\cdots +n_k+n^*_{k+1}+\cdots +n^*_{k+j}$ for some $j\in\{1\ld l\}$. Then $w(\ep_i)=\ep_{i+1}$ unless $i=r(j)$ for some $j$, and $w(\ep_{r(j)})=-\ep_{r(j-1)+1}$.

Let $\GC=SO_{2n}(\overline{F}_q)$, $q$ odd, or $Spin_{2n}(\overline{F}_q)$, $q$ even. Then the conjugacy classes of $W$ are parameterized by
 double partitions  $[n_1\ld n_k,n^{*}_{k+1}\ld n^{*}_{k+l}]$ with  $l$ is even, except for the cases where $l=0$ and each number $n_1\ld n_k$ is even \cite[Prop 25]{Ca1}. In the exceptional cases there are two conjugacy classes
 corresponding to $[n_1\ld n_k,-]$. For a  canonical representative of one class we choose the permutation $w$
 of $\ep_1\ld \ep_n$ described above for  $GL_{n}(\overline{F}_q)$.  A  canonical representative  $w'$ of the second class is chosen  as follows. We  set  $w'(\ep_1)= -w(\ep_1)$, $w'(\ep_{n_1})= -w(\ep_{n_1})$   and $w'(\ep_i)= w(\ep_i)$ for $i\neq 1,n_1$ ($1\leq i\leq n$).

\subsection{Canonical tori}

(A)
Suppose first that $G=GL_n(q)$ and $\TC_1$ is the group of diagonal matrices. Then $F(t)=t^q$ for $t\in \TC_1$. Let $w\in W\cong S_n$ be  a canonical  representative of   the conjugacy class corresponding  to a partition $[n_1\ld n_k]$ of $n$. If $k=1,n_1=n$, then $w(i)=i-1$ for $i=2\ld n$, and $w(1)=n$.
 Then $F(t)=t^q=w(t)$ \ii $t=(b,b^q\ld b^{q^{n-1}})$ for $b\in   \overline{F}_q$ with $b^{q^{n}-1}=1.$ This is exactly the canonical $w$-torus $D_w$ for this $w$. 
 If one fixes a primitive
$(q^{n}-1)$-root of unity $a$,  then $D_w$ is the group generated by  $(a,a^q\ld a^{q^n-1})\in \TC_1$. 
This special case illustrates for a reader what we do in general. 

For an arbitrary partition $\pi=[n_1\ld n_k]$ we consider the group of diagonal matrices $\diag(D_1\ld D_k)$, where $D_i$ is generated by a matrix $\diag(a_i, a_i^q\ld a_i^{q^{n_i-1}})$ and $a_i$ is an arbitrary primitive $(q^{n_i}-1)$-root of unity. 
One observes that this group coincides with $D_w$ for the canonical choice of $w=w(\pi)$ as described above.

\medskip
(B) The situation with other classical groups is similar but requires adjustments. We start with groups
$\GC=SO_{2n+1}(\overline{F}_q)$ with $q$ odd, or $Sp_{2n}(\overline{F}_q)$. Then the conjugacy classes of $W$
are in bijection with the double partitions $\pi=[n_1\ld n_k,n_{k+1}^*\ld n^*_{k+l}]$.

Suppose first that     $k=0,l=1$. For $t=(x_1\ld x_n)\in \TC_1$
we have $w(t)=(x_2\ld x_n,x_1\up)$ and $F(t)=(x^q_1\ld x^q_n)$. The equality $w(t)=F(t)$ implies
$x_i=x_{i-1}^q$ for $i=2\ld n$ and $x_1\up=x_n^q$. So $x_n=x_1^{q^n}=x_1\up$, whence $x_1^{q^n+1}=1$.
\itf  the set $\{t\in \TC_1:w(t)=F(t)\}$ coincides with $\{(b,b^q\ld b^{q^{n_i}}): b^{q^n+1}=1\}$. In particular,  this is a cyclic group of order $q^n+1$.

In general,  the group $D_w=\{t\in \TC_1:w(t)=F(t)\}$ is the direct product of groups $D_1\ld D_{k+l}$, where
$D_i$ is a cyclic group of order $q^{n_i}-1$ for $i\leq k$ and of order $q^{n_i}+1$ for $i> k$. 
Specifically, if $r=n_1+\cdots +n_{i-1}$ then $D_i=(1\ld 1,b_i,b_i^q\ld b_i^{q^{n_i-1}},1\ld 1)$, 
where $b_i$ occupies $(r+1)$-th position and $b_i^{q^{n_i}-1}=1$ for $i\leq k$ and  
$b_i^{q^{n_i}+1}=1$ for $i>k$.

(C)  Let $\GC=SO_{2n}(\overline{F}_q)$, $q$ odd, or  $Spin_{2n}(\overline{F}_q)$, $q$ even , $W=W_\GC$ and for a moment let $\tilde W$ be the Weyl group of $SO_{2n+1}(\overline{F}_q)$.    
 Then $W$ is a normal subgroup of $\tilde W$, and  an element of  $\tilde W$  lies in $W$ \ii the conjugacy class of it  corresponds to the double partition $\pi$
with $l$ even  \cite[Prop 25]{Ca1}. In particular, the canonical element from such a class lies in $\tilde W$. \itf  
 %If $l$ is even then the canonical element $w$ corresponding to $\pi$ belongs to the Weyl group of $SO_{2n}(\overline{F}_q)$. 
  $D_w$  is the same for groups $SO_{2n}(\overline{F}_q)$ and $SO_{2n+1}(\overline{F}_q)$ for $w\in W$, as well as for the pair $Spin_{2n}(\overline{F}_q)$, $q$ even, and  $Sp_{2n}(\overline{F}_q)$, $q$ even.  
 However in the exceptional cases, when $l=0$ and all $n_1\ld n_k$ are even, the Weyl group of the former group 
 has two conjugacy classes corresponding to $\pi=[n_1\ld n_k,-]$, so we need to construct a canonical $w$-torus   for $w$ from the second   class corresponding to   $\pi$.

Consider first a special case with $\pi=[n,-]$, $n$  even. If $w$ permutes $\{1\ld n\}$ 
then the group $\{(b,b^q\ld b^{q^{n-1}}): b\in F_q,b^{q^{n}-1}=1 \}$ satisfies $w(t)=F(t)=t^q$
if $w(1)=n$, $w(i)=i-1$ for $i<n$. 
 A canonical  representative $w'$  of the other conjugacy class 
corresponding to $\pi$  is defined to be 
$w'(\ep_i)=w(\ep_i)$ for $i\neq 1,n_1$, and   $w'(\ep_{1})=-w(\ep_{1})$,  $w'(\ep_{n_1})=-w(\ep_{n_1}).$
So
$w(x_1\ld x_n )=(x_2\up,x_3\ld x_{n-1},x_1\up)$.  Let  $t=(b\up,b^{q},b^{q^2}\ld b^{q^{n-1}})\in \TC_1$ with $b\in F_q$ and $b^{q^n-1}=1$. Then 
$w'(t)=F(t)=t^q$, which  implies that
$t=(b\up,b^{q},b^{q^2}\ld b^{q^{n-1}})$ for some $b\in F_q$ with $b^{q^n-1}=1$.

Let $k>1$ and $w$ correspond to the double partition $\pi=[n_1\ld n_k,-]$ with
 even parts $n_1\ld n_k$. If $w$ just permutes $\ep_1\ld \ep_n$ then  $D_w$ is as above, that is,
 the direct product of subgroups $D_1\ld D_k$, where
 $D_i=(1\ld 1,b_i,b_i^q\ld b_i^{q^{n_i-1}},1\ld 1)$ ($1\leq i \leq k$), 
 $b_i$ occupies the $(n_1+\cdots +n_{i-1}+1)$-th position and $b_i^{q^{n_i}-1}=1$. If $w'$ is a canonical representative of the other 
conjugacy class  of $W$ corresponding to $\pi$ then $D_{w'}$
is the direct product of subgroups $D'_1, D_2\ld D_k$, where
$D_2\ld D_k$ are as above, and $D'_1=\{b_1\up,b_1^q\ld b_1^{q^{n_1-1}}\}$
with $b_1^{q^{n_1}-1}=1.$

(The groups $D_w,D_{w'}$ are known to be conjugate in
 $O(2n,\overline{F}_q)$  but not in  $SO_{2n}(\overline{F}_q)$.)

  \medskip
 (D)  It is well known that $G^+=SO^+_{2n}(q)$, $q$ odd,  can be viewed as a subgroup of $H=SO_{2n+1}(q)$, which  agrees
  with an  inclusion $\GC\subset {\mathbf H}=SO_{2n+1}(\overline{F}_q)$, in the sense that the Frobenius morphism defining $G^+$ is the restriction to $\GC$ of a Frobenius morphism $F$ of ${\mathbf H}$ defining $H$. (Similarly, for  the pair $Spin_{2n}(\overline{F}_q)$, $q$ even, and  $Sp_{2n}(\overline{F}_q)$.) Then the maximal reference torus $\TC_1$ of ${\mathbf H}$
coincides with $\TC_1$ for  $\GC$, and $F(t)=t^q$ for $t\in \TC_1$ for both the groups. This  yields an embedding $W_\GC=N_{{\mathbf G}}(\TC_1)/\TC_1$ into   
     $W_{\mathbf H}=N_{{\mathbf H}}(\TC_1)/\TC_1$.
  
  Denote for a moment by $F^+$ and $F^-$ the Frobenius morphisms of  
  $\GC$ such that $\GC^{F^+}=G^+=SO^+_{2n}(q)$ and $\GC^{F^-}=G^-=SO^-_{2n}(q)$. Then $F^-=\si \cdot F^+$, where $\si$ is a graph \au of $\GC$ \cite[\S 11]{St}.    Moreover, $\si^2=1$ and $\si$ can be realized
  as a conjugation by an element $d$ of $N_{{\mathbf H}}(\TC_1)$, that is, if $x\in\GC$ then $\si F^+(x)=d(F^+(x))d\up$. Let $r$ be the projection of $d$ into $W_{\mathbf H}$.  
As $F$, and hence $F^+$, act trivially on $W_{\mathbf H}$, the action of $ F^-$ on 
$W_{\mathbf G}$ coincides with the conjugation by  $r$. Clearly, $r\notin W_{\mathbf G}$. As $W_{\mathbf G}$ has index 2 in $W_{\mathbf H}$, we have $W_{\mathbf H}=W_{\mathbf G}\cup rW_{\mathbf G}$, and the coset $rW_{\mathbf G}$ is invariant in $W_{\mathbf H}$.   

We need to describe canonical $w$-tori $D_w$  for $G^-$. For our purpose it is  convenient to 
  describe them in terms of group $H$.   To avoid confusion, we use notation $D^-_w$ for canonicall $w$-tori of $G^-$ and keep $D_w$ for those  of $H$. 
 Note that, for uniformity reason,  in Lemma \ref{nn6} and Proposition \ref{pp1} we allow $q$ to be even when dealing with the group $SO_{2n+1}(q)\cong Sp_{2n}(q)$.
  
  \bl{nn6} $D^-_w=D_{rw}$, and we use for $D^-_w$ the canonical $rw$-torus $D_{rw}$ 
  of $H$. 
  \el
  
  Proof. We have $D^-_w=\{t\in \TC_1: F^-(t)=w(t)\}$. As $F^-(t)=rF^+(t)=rF(t)$ and $r^2=1$,
 we have  $F^-(t)=rF^+(t)=rF(t), $ so the equality $F^-(t)=w(t)$ is equivalent to $F(t)=rw(t)$, whence the result. 
  
   \medskip
   We summarize the considerations of this section as follows:

\begin{propo}\label{pp1}  $(1)$ Let $\GC^F=G=Sp_{2n}(q)$ or $SO_{2n+1}(q)$. Let $w\in W_\GC$ be a canonical element corresponding to the double partition $\pi=[n_1\ld n_k,n^*_{k+1}\ld n^*_{k+l}]$. Then    $\TC_1^{w^{-1}\circ F}=D_w$,  where $D_w=(D_1\ld D_{k+l})$,
 as described above.

$(2)$ Let  $\GC^F=G=SO^+_{2n}(q)$  and $w\in W_\GC$ be a canonical element corresponding to $\pi$. 
If $\pi$      is non-exceptional then $D_w$ is as above. If $\pi$      is  exceptional  then  
there is one more canonical element $w'$ corresponding to $\pi$ 
that is not $W_\GC$-conjugate to w. Then $D_w$ is   as above, and   
$D_{w'}=(D'_1\ld D_k )$, where $D_1'=(b_1\up,b_1^q\ld b_1^{q^{n_1-1}})$ with $b_1\in \overline{F}_q$ and $b_1^{q^{n_1}-1}=1$. 

$(3)$  Let  $G=SO^-_{2n}(q)$. Then $D^-_w$ coincides with the canonical torus $D_{v}$ of ${\mathbf H}=SO_{2n+1}(\overline{ F}_q)$ for some  $v\in W_{\mathbf H}\setminus W_\GC$.
\end{propo}

\subsection{Weights and $q$-characters} %maximal tori}

Let $\GC$  be an algebraic group with Frobenius morphism $F$, and $G=G^F$. 
For certain weights $\mu$ and every maximal torus $T$ of $G$ we compute
$(s_\mu|_T,1_T)$, the multiplicity of the trivial character of $T$ in the orbit character $s_\mu$.
It is well known that $(s_\mu|_T,1_T)$ is the same for any $\GC$-conjugate of $T$ in $\GC$.
Therefore, it suffices to do this for the canonical representative $D$ of $T$ in $\TC_1 $.

Below $\GC\in\{GL_n(\overline{F}_q), SL_n(\overline{F}_q), $ $Sp_{2n}(\overline{F}_q)$,  $SO_{2n+1}(\overline{F}_q)$, $q$ odd,     $\GC=SO_{2n}(\overline{F}_q)$, $q$ odd,  $Spin_{2n}(\overline{F}_q)$,
$q$ even$\}$. The weights $\ep_1\ld \ep_n$ and $\om_1\ld \om_n$ are defined in the introduction.  So
$G=\GC^F$ is a classical group (except  a unitary one, which is not  considered).  As above, $W$ is the Weyl group of $\GC$ and $W_j=C_{W}(\om_j)$.  We write $W\om_j$ for the $W$-orbit of $\om_j$.

\begin{defi}\label{d91} %\mar{d91}
 Let $T$ be a maximal torus of $G$, and $\theta\in\Irr T$. We say that $\theta$ is a $q$-character if $\theta(t)^{q-1}=1$ for all $t\in T$. This also applies to Brauer characters $T\ra \overline{F}_q$. 
\end{defi}

\bl{t11} %\mar{t11}
Let $\mu $ be a weight of $\GC$ and T a maximal torus of G. Then $\mu|_T$ is a q-character of T \ii $((q-1)\mu)|_T=1_T$. \el

Proof. If $\mu$ is a $q$-character of $T$ then, obviously,   $((q-1)\mu)|_T=1_T$.
 Conversely, suppose that $((q-1)\mu)|_T=1_T$. Let $t\in T$, so a weight $\mu$ is  a $q$-character of $T$ \ii   $((q-1)\mu(t))=1$, equivalently, $\mu(t)^{q-1}=1$. As $\mu(t)\in \overline{F}_q$, and $x^{q-1}=1$ for $x
\in \overline{F}_q$ implies $x\in F_q$, we have $\mu(t)\in F_q$, whence the claim.

\begin{lemma}\label{au1}\mar{au1} Let $1\leq k\leq n$, $0<r<q$ and
$0\leq l_1<\cdots <l_k<n$ be integers. Then $r(q-1)(q^{l_1}+\cdots
+q^{l_k})$ is not divisible by $q^n-1$, unless $k=n$ and $(l_1\ld
l_{n})=(0,1\ld n-1) $.\end{lemma}

Proof. Suppose first that $l_k\leq n-2$, so $k\leq n-1$.
Then $r(q-1)(q^{l_1}+\cdots +q^{l_k})\leq r(q-1)(q^{n-k-1}+\cdots
+q^{n-3}+q^{n-2})=r(q^{n-1}-q^{n-k-1})$. As $r<q$, we have
$r(q^{n-1}-q^{n-k-1})<q^n-q^{n-k}=q^n-1 -(q^{n-k}-1)$, which is
less than $q^n-1$. So we are left with the case   $l_k=n-1$.

Let $l_{c+1}\ld l_{k-1},l_k=n-1$ be the longest string of
subsequent natural numbers, that is,  $(l_{c+1}\ld l_k)=(n-(k-c)\ld n-1)$ and $l_c<n-(k-c)-1$.
So  $c+1 \leq k$. 

If $c=0 $ then $(l_1\ld
l_{k})=(n-k\ld n-1)$  and $r(q-1)(q^{l_1}+\cdots
+q^{l_k})=r(q^n-q^{n-k})$, so the lemma follows if $q^{n-k}>1$.
If $q^{n-k}=1$ then $n=k$ and $(l_1\ld l_{k})=(0\ld n-1)$, so we are in
the exceptional case of the lemma.)  So we assume $c>0 $.

 Then
$l_j\leq n-2 -(k-j)$ for $1\leq j\leq c$. Therefore,
$(q-1)(q^{l_1}+\cdots +q^{l_c})\leq (q-1)(
q^{n-k-1}+q^{n-k-2}+\cdots +q^{n-2 -(k-c)})= q^{n-1 -
(k-c)}-q^{n-k-1}<q^{n-1-(k-c)}$ as $c>0$. We have
$$r(q-1)(q^{l_1}+\cdots +q^{l_k}) = %r(q-1)(q^{l_1}+\cdots
%+q^{l_{c}}+q^{l_{c+1}}+\cdots +q^{l_k})=
r(q-1)(q^{l_1}+\cdots +q^{l_{c}}+q^{n-(k-c)}+\cdots +q^{n-1})=$$
$$= r(q-1)(q^{l_1}+\cdots +q^{l_{c}})+r(q^n-1)-r(q^{n-(k-c)}-1).$$
This is a multiple of $q^n-1$ \ii so is $x:=r(q-1)(q^{l_1}+\cdots
+q^{l_{c}})-r(q^{n-(k-c)}-1)$. It is easy to observe that $x$ is not
divisible by $q^n-1$. (Use the inequality $r(q-1)(q^{l_1}+\cdots
+q^{l_{c}})\leq r(q^{n-1 -(k-c)}-q^{n-k-1})<r(q^{n-(k-c)}-1)$.)

\medskip
(A) $\GC=GL_{n}(\overline{F}_q)$ or $SL_{n}(\overline{F}_q)$. 

Let  $\pi=[n_1\ld n_k]$ be a partition of $n$ and $w\in W$ the canonical element in the conjugacy class determined by $\pi$. Let $B_1\ld B_k$ be the orbits of $w$ on $\{1\ld n\}$. Then $\{1\ld n\}=B_1\cup \cdots \cup B_k$, where $B_1=\{b: 0<b\leq n_1\}$ and $B_i=\{b: n_1+\cdots +n_{i-1}<b\leq n_1+\cdots +n_{i-1}\}$ for $i=2\ld k$. 

Let   $D_w=(D_1\ld  D_k)$ be a canonical $w$-torus in $\TC_1$. Recall that   $D_i$ $(1\leq i\leq k)$ is generated by an element $ (a_i,a_i^q,\ld a_i^{q^{n_i-1}})$, where  $a_i\in \overline{F}_q$
is a primitive $(q^{n_i}-1)$-root of unity.  In this notation we have: %We consider the characters of $D_w$ of shape $(\ep_{m_1}+\cdots+ \ep_{m_j})|_{D_w}$. 

\begin{propo}\label{n1n} %\mar{n1n}
 Let  $G=GL_n(q)$ or  $SL_n(q)$.  
Set   $\mu=g(\om_j)$ for some $g\in W$.

$(1)$ $\mu|_{D_w}$ is a $q$-character \ii $w(\mu)=\mu;$

$(2)$ %the W-orbit of $\om_j$. 
The number of distinct q-characters $\mu|_{D_w}$  $(\mu\in W\om_j)$ equals $1_{W_j}^W(w)$. In addition, $1_{W_j}^W(w)=(s_{(q-1)\om_j}|_T,1_T)$.\end{propo}

Proof.  
(1) The orbit $W\ep_j$ consists of all weights $\ep_{m_1}+\cdots +\ep_{m_j}$ $(1\leq m_1<\cdots <m_j\leq n)$. Let $\mu=\ep_{m_1}+\cdots +\ep_{m_j}.$

Suppose first that $k=1$, so $B=B_1=\{1\ld n\}$. If   $G=GL_n(q)$ then  $D_w=\lan d\ran$ is a cyclic group of order $q^n-1$ with generator  
 $d= (a,a^q\ld a^{q^{n-1}})$,
where $a\in \overline{F}_q$ is a primitive $(q^n-1)$-root of unity.  
Then $\mu(d)=a^{q^{m_1-1}+\cdots +q^{m_j-1}}$. This belongs to $F_q$
\ii $(q-1)(q^{m_1-1}+\cdots
+q^{m_j-1})\equiv 0\pmod {(q^n-1)}$.  If   $G=SL_n(q)$ then $D_w=\lan d^{q-1}\ran$ is of order  $(q^{n}-1)/(q-1)$. So $(\ep_{m_1}+\cdots +\ep_{m_j})(d^{q-1})\in F_q$ \ii 
 $(q-1)(\ep_{m_1}+\cdots +\ep_{m_j})(d)\in F_q$,
 equivalently, $(q-1)^2(q^{m_1-1}+\cdots +q^{m_j-1})\equiv 0\,(q^n-1)$.  By Lemma \ref{au1}, applied to
 $\{m_1-1\ld m_j-1\}$ in place of $\{l_1\ld l_k\}$ with $r=1$ or $q-1$, this is not the case unless $j=n$ and $(m_1\ld
m_{j})=(1\ld n) $.

Let $k>1$.  Set $B_i'=\{b\in B_i: b\in \{m_1\ld m_j\}\}$. 
 Let $G=GL_n(q)$ and let $d_i$ be a generator of
 the subgroup $D_i=(\Id,\ldots,\Id,D_i,\Id\ld \Id )$ of $D_w$. Then  $\mu(d_i)=\sum_{r\in B_i'} \ep_{r}(d_i)$. The above argument shows that 
 this is in $F_q$ \ii $B_i'=B_i$. As this is true for every $i$ with non-empty $B_i'$, it follows that
 $\{m_1\ld m_j\}$ is the union of $w$-orbits, whence (1).

Let $G=SL_n(q)$. Then  $D_w$ contains a subgroup which is the direct product of cyclic subgroups $D_i^{q-1}$ ($1\leq i\leq k$), where $D_i=\lan d_i\ran$ is as above for $GL_n(q)$.    Suppose that $\mu(D_w)=1$.    Then $\mu(d_i^{q-1})=\sum _{r\in B_i'}\ep_{r}(d_i^{q-1})=1$. By the argument for $k=1$, this implies  
$B_i'=B_i$, so again $\{m_1\ld m_j\}$ is the union of $w$-orbits, as desired.

(2) Note that $g(\om_j)=h(\om_j)$ ($g,h\in W$) \ii $gW_j=hW_j$. Let $\mu=g(\om_j)$. 
Then $w(\mu)=\mu$, or $wg(\om_j)=g(\om_j)$, is equivalent to $wgW_j=gW_j$. So the number of distinct $q$-characters   $g(\om_j)$ is equal to the number of $w$-stable cosets $gW_j$. This is well known to be equal to the value at $w$ of the character $1_{W_j}^W$, as claimed.

 \medskip
 Remark. Strictly speaking, the proof of Proposition \ref{n1n} does not require the choice of $w$
 to be canonical, which only affects  the explicit form of the sets $B_1\ld B_k$. 
  
\begin{corol} Let $G=SL_n(2)$ and let $V_j$ $(j<n)$ be the j-th exterior power of the natural $F_2G$-module V.  Let $T=T_w$ be a maximal torus of G, and $V^{T_w}_j$ the fixed point subspace  
of $T_w$ on $V_j$. Then $\dim V^{T_w}_j =1_Y^{S_n}(w)$, where   $Y$ is the Young subgroup of $S_n$ labeled by $[j,n-j]$.
\end{corol}

Proof. Let $\rho$ be the \irr of $\GC =SL_n(\overline{F}_2)$ with highest weight $\om_j=\ep_1+\cdots +\ep_{n-1} $. As $\om_j$ is a miniscule weight \cite[\S 7.3]{Bo2}, the weights of $\rho$
are  $g(\om_j)$ for $g\in W$.  As $\dim V^{T_w}_j $ equals the number of  distinct weights $g(\om_j)$
such that $g(\om_j)(T_w)=1$, we have $\dim V^{T_w}_j =1_{W_j}^{W}(w)$. Since $W\cong S_n$ and $W_j\cong Y$, the statement follows.

\medskip
(B) Let $G\in\{Sp_{2n}(q), SO^+_{2n}(q), ~q~odd,~ SO_{2n+1}(q),~q~odd,~Spin_{2n}^+(q),~q~even \}$, 
%SO^-_{2n}(q)$ 
and $\GC$ the respective algebraic group.
  We use Proposition  \ref{pp1} to compute, for every  canonical $w$-tori  $D_w$, the number of 
weights $\mu$ in the orbit $W\om_j$ whose restriction to $D_w$ yields   a $q$-character of $D_w$. 
  It is well known that  $W\om_j$ consists of all
weights $\pm \ep_1\pm \cdots \pm \ep_{m_j}$, except if $\GC$ is of type $D_n$  
and $j=n$, where the orbit
consists of all weights $\pm \ep_{1}\pm \cdots \pm \ep_{n}$ with even number of minus signs.
So the $W_\GC$-orbit of $\om_j$ is the same for 
the above groups $\GC$, except for the case where $j=n$ and $\GC$ is  of type $D_n$.

 Below the exceptional  $w$-torus  $D_w$ is one in 
 $SO^+_{2n}(q)$ for the exceptional canonical element $w\in W$; the corresponding partition $\pi$  is $[n_1\ld n_k,-]$, where all $n_1\ld n_k$ are even. Observe that there are one exceptional and one non-exceptional canonical elements corresponding to this partition but the exceptional canonical $w$-torus exists only in $SO^+_{2n}(q)$. Sometimes we denote the exceptional canonical element by $w'$ and the corresponding torus by $D_{w'}$ to avoid confusion. 

\bl{q-sp}\mar{q-sp} Let $G\in\{Sp_{2n}(q),$ q odd, $ SO_{2n+1}(q), SO^+_{2n}(q), q~odd,    Spin^+_{2n}(q)$, $q$ even$\}$, and 
 $D_w$  the canonical $w$-torus for $G$ corresponding to  
 $w\in W$.  Let $\mu\in W\om_j$. Then the \f conditions are equivalent:

$(1)$ $\mu|_{D_w}$ is a q-character of $D_w;$

$(2)$ $w(\mu)=\mu$.

%\noindent In addition, if $j=n$ and $G=SO^-_{2n}(q)$ then $\mu|_{T}$ is not a q-character for every maximal torus T of $G$.
\el

Proof. Let $\pi=[n_1\ld n_k,n^*_{k+1}\ld n^*_{k+l}]$ be a double partition corresponding to $w$. Denote by $B_i$ ($1\leq i\leq k+l$) the set of natural numbers
$r$ in the range $n_1+\cdots +n_{i-1}<r\leq n_1+\cdots +n_i$. So $|B_i|=n_i$.
We use the canonical form of $D_w$ described in Proposition  \ref{pp1}. So $D_w=(D_1\ld D_{k+l})$ and in the exceptional case (where $G=SO^+_{2n}(q)$, $q$ odd, or $Spin^+_{2n}(q)$, $q$ even, $l=0$ and all  $n_1\ld n_k$ are even) we also
consider $D_{w'}=(D'_1,D_2\ld D_{k})$. Note that $D_i$ is generated by an element  $d_i=(1\ld 1,a_i,a_i^q$ $\ld a_i^{q^{n_i-1}},1\ld 1)$ and $D_1'$ is generated by $d_1'=(a_1\up,a_1^q\ld a_1^{q^{n_1-1}},1\ld 1)$.
%Let $\pi$ be a partition such that $w$ is on conjugate class $C_\pi$ of $S_n$.

$(i)$ Suppose that $\pi=[n,-]$ or $[-, n^*]$.  % $k=1,l=0$, so %$B=B_1=\{b_1\ld b_n, b_1'\ld b_n'\}$  and 
Then  $D_w=\lan d\ran$ is a cyclic group of order $q^n-1$ or  $q^n+1$, respectively, and 
 $d:= (a,a^q,a^{q^2}\ld a^{q^{n-1}})\in \TC_1$ for $a\in \overline{F}_q$
with $|a|=|d|$. If $\pi=[n,-]$ with $n$ even and $G=SO^+_{2n}(q)$, $q$ odd, or  $Spin^+_{2n}(q)$, $q$ even,  then there is also 
an exceptional canonical element $w'$ for which 
$D_{w'}=\lan d'\ran$, where
 $d':= (a\up ,a^q,a^{q^2}\ld a^{q^{n-1}})$.

  Recall that $\mu=\pm\ep_{m_1}\pm\cdots\pm \ep_{m_j}$ $(1\leq m_1<\cdots <m_j\leq n) $  for a  certain choice of signs.
So  
$\mu(d) =a^{\pm q^{m_1-1}\pm \cdots \pm q^{m_j-1}}$, as well as $\mu(d')$ 
 in the exceptional case. %with $m_1=1$ we have $(\pm\ep_{m_1}\pm\cdots\pm \ep_{m_j})(t)=a^{-\pm q^{m_1-1}\pm \cdots \pm q^{m_j-1}}$. 
 Note that  $a^{\pm q^{m_1-1}\pm \cdots \pm q^{m_j-1}}\in F_q$ \ii
$(q-1)(\pm q^{m_1-1}\pm \cdots \pm q^{m_j-1})\equiv 0 \pmod {|a|}$. The  left hand side  is not 0 as   $m_1\ld m_j$ are distinct. Furthermore, the absolute value of the  left hand side 
is  strictly less than $|a|$,  unless 
$\{m_1\ld m_j\}=\{1\ld n\}$, $|a|=q^n-1$ and the signs of all $\pm 1, \pm q,  \ldots ,\pm q^{n-1}$ are the same. 
In particular,  $\mu(d)\notin F_q$ if $j<n$ or $j=n$ and $|a|=q^n+1$. %This justifies the additional statement for our  choice of $w$. 

Let $j=n$, $|a|=q^n-1$. If we compute $\mu(d)$ (respectively, $\mu(d')$) then the signs in  $\pm 1  \ldots ,\pm q^{n-1}$ are the same 
\ii  $\mu=\pm(\ep_1+\cdots  +\ep_n)$ (respectively, 
$\mu=\pm(-\ep_1+\ep_2+\cdots  +\ep_n))$. 

If  $\GC=SO_{2n}(\overline{F}_q)$, $q$ odd, or $Spin_{2n}(\overline{F}_q)$, $q$ even, then  the weights $\pm(-\ep_1+\ep_2+\cdots  +\ep_n)$ are not in $W\om_n$, and this possibility is ruled out. However, we use below the above statement for $\mu(d')$.

Thus, if (1) holds then $j=n$
and $w$ permutes $\ep_1\ld \ep_n$, so (2) holds. 

Conversely, if $w(\mu)=\mu$ then $w$ is not of type $[-,n^*]$ as otherwise $w^n(\mu)=-\mu$.
So $w$ is of type $[n,-]$, and $|a|=q^n-1$, which implies (1) with the above observations.

$(ii)$ The general case. 

$(1)\ra (2)$ Suppose  that %$w$ is not exceptional. Then  
$\mu(D_w)\subset  F_q$. Then $\mu(d_i)\in F_q$ for every $i=1\ld k+l$
in the non-exceptional case, otherwise  $\mu(d_1')\in F_q$ and  $\mu(d_i)\in F_q$ for every $i=2\ld k$.  Observe that  $\mu(d_i)=\sum_r \pm \ep_r(d_i)$, where $r$ runs over $B_i':=\{m_1\ld m_j\}\cap B_i$, and similarly for $d_1'$.   Suppose that $B_i'$ is not empty.    Then the same reasoning as above with $a_i$ in place of $a$ and $\sum_{j\in B_i} \ep_j$ in place of $\ep_1+\cdots +\ep_n$ shows that  $i\leq k$ and $B_i'=B_i$ (so $\{m_1\ld m_j\}\subseteq B_1\cup\cdots\cup B_k$).

Furthermore, the argument in (i) tells us that the signs of $\ep_{r}$ are the same for 
all  $r\in B'_i$ in the non-exceptional case. In the exceptional case this is only true for $i>1$, 
whereas for $i=1$ all but one signs of $\ep_r$ with $r\in B_1'$ are the same.

    In the   non-exceptional case  $w$ simply permutes $\ep_r$ with $r\in B_i$ for $i\leq k$, so $w(\mu)=\mu$ as required. 

Consider  the exceptional canonical torus $D_{w'}=(D_1',D_2\ld D_k)$ for $G=SO^+_{2n}(q)$, $q$ odd, and $Spin^+_{2n}(q)$, $q$  even, so $l=0$, $|B_i|=n_i$ are even and $|d_i|=q^{n_i }-1$.  
 Then  $d_1'=(a_1\up,a_1^q\ld a_1^{q^{n_1-1}},1\ld 1)$, whereas  $D_i=\lan d_i\ran$ with $i>1$.
%is generated by $d_i=(1\ld 1,a_i,a_i^q$ $\ld a_i^{q^{n_i-1}},1\ld 1)$. 
If $B'_1$ is empty 
then $\{m_1\ld m_j\}$ is the union of $B_i$ with $1<i\leq k$ and the signs of $\ep_{r}$ are the same for 
all  $r\in B'_i$. Suppose that $B_1'$ is non-empty. 
Then $B_1\subseteq \{m_1\ld m_j\}$, so $n_1\leq j$ and  $(m_1\ld m_{n_1})=(1\ld n_1)$. 
 As above, if $\mu(d'_1)\in F_q$ then the partial sum
of 
$\mu=\pm \ep_{1}\pm \cdots \pm \ep_{n}$ with indices in $  B_1$ must be $\pm (-\ep_1+\ep_2+\cdots +\ep_{n_1})$. In addition, the partial sums corresponding to  $B_i$ with $i>1$ %, that meet $\{m_1\ld m_j\}$ 
must be $\pm (\sum_{r\in B_i}\ep_r)$. If $j<n$ then  %there exists $g\in W$ such that
%By the above (subsect 3.2?) 
$w'(\ep_1)=-\ep_2$, $w'(\ep_{n_1})=-\ep_1$, and $w'(\ep_i)=\ep_{i+1}$ for $i=2\ld n_1-1$,
and also $w'(\sum_{r\in B_i'}\ep_r)$ for $i>1 $.
Therefore, $w'(\pm (-\ep_1+\ep_2+\cdots +\ep_{n_1}))=\pm (-\ep_1+\ep_2+\cdots +\ep_{n_1})$, and (2) follows.

Let $j=n$.  Recall that $|B_i|=n_i$ is even for every $i$.   Therefore, the number of minus signs  in the expression  $\mu=\pm \ep_{1}\pm \cdots \pm \ep_{n}$ is odd. However, such expression cannot be in the $W$-orbit of $\ep_1+\cdots +\ep_n$, which is a contradiction. Therefore, this case does not occur, and hence $(1)\ra (2)$ for arbitrary $w\in W$.

$(2)\ra (1)$   Let $w(\mu)=\mu$. Then, obviously, the $w$-orbit of each $\ep_{m_r}$ $(r=1\ld j)$ is in the set $\{\pm \ep_{m_1}\ld \pm \ep_{m_j}\}$, so $B_i'=B_i\cap \{m_1\ld m_j\}\neq \emptyset $ implies $B_i\subseteq \{m_1\ld m_j\}$.
Then $i\leq k$ (as otherwise for a partial sum $\nu:=\sum_{r\in B_i} \pm \ep_{m_r}$ we have $w^{n_i}(\nu)=-\nu$, and hence $w^{n_i}(\mu)\neq \mu$). If $w$ is non-exceptional, this  
 means that $w$  permutes $\ep_r$ for $r\in B'_i$ for every $i$ (with $B_i'\neq\emptyset$), so the signs of 
 these $\ep_r$'s  are the same. This implies that $\mu(D_w)\subset F_q$. Suppose that $w$ is exceptional.
This argument again works if $B_1'=\emptyset$. Otherwise,  we have to pay an additional attention to 
the case where $w=w'$ and $D_{w'}=(a_1\up, a_1^q\ld a_1^{q^{n_1}-1}, a_2\ld a_2^{q^{n_2-1}}\ld a_k\ld a_{k}^{q^{n^k-1}})$.  Then $w'(\mu)=\mu$ implies $\mu=\pm (-\ep_1+\ep_2+\cdots +\ep_{n_1})+\sum_{i=2}^k\pm(\sum_{r\in B_i'}\ep_r)$. If $j=n$ then this $\mu$ is not in the $W$-orbit of $\om_n$.
If $j<n$ then $\mu(D_w)\subset F_q$, as required.

\bl{hb4} %\mar{hb4}
 $(1)$ Let $G\neq SO^-_{2n}(q)$, q odd,  $Spin^-_{2n}(q)$, $q$ even. 
  For $w\in W$ let $T_w$  be a maximal torus of $G$.  Let $W_j=C_{W}(\om_j)$. 
Then the number of 
distinct weights $g(\om_j)$  $(g\in W)$ that yield     q-characters of $T_w$  equals  the number of the cosets $gW_j$  such that 
$wgW_j=gW_j$. The latter equals $1_{W_j}^{W}(w)$.

 $(2)$  Let $G=SO^-_{2n}(q)$, q odd,     $H=SO_{2n+1}(q)$, $q$ odd, $\tilde W=W_{{\mathbf H}}$
 and $\tilde W_j=C_{\tilde W}(\om_j)$ for $j<n$.  Let $r\in \tilde W$ be as in Lemma \ref{nn6}. 
 Then the  number of distinct weights  
 $ g(\om_j)$ $(g\in W)$ that yield q-characters of $T_w$ equals $1_{\tilde W_j}^{\tilde W}(wr)$, where r is as in Lemma {\rm \ref{nn6}}. 
  (Similarly, for  $G=Spin^-_{2n}(q)$, $Spin_{2n}^-(q)$, for $q$ even.)
\el

Proof. It suffices to proof the lemma for the canonical $w$-torus $D_w$ in place of $T_w$. (1) Let $\mu=g(\om_j)$.  By Lemma \ref{q-sp},  $\mu|_{D_w}$ is a $q$-character \ii $w(\mu)=\mu$, equivalently, $(g\up wg)(\om_j)=\om_j
$, that is, $g\up wg\in W_j$. Let $h\in W$.
%Recall that the $W$-orbit of $\om_j$ consists of all weights $\pm\ep_{m_1}\pm\cdots\pm \ep_{m_j}$, except for the case where $j=n$ and $\GC=SO(2n, \overline{F}_q)$.  Set $\mu=\pm\ep_{m_1}\pm\cdots\pm \ep_{m_j}$. By Lemma \ref{q-sp},  $\mu|_{T_w}$ is a $q$-character \ii $w(\mu)=\mu$. Let $g\in W$ be such that $g(\om_j)=\mu$. 
Then $h(\om_j)=\mu$  \ii  $h\in gW_j$. \itf  the number of 
$q$-characters in the set $(g(\mu))|_{D_w}$ is equal to the number of 
cosets $gW_j$  such that $wgW_j=gW_j$. This is exactly $1_{W_j}^W(w)$.   

(2) By Lemma \ref{nn6}, $D_w=\tilde D_{wr}$, where $\tilde D_{wr}$ is the canonical $wr$-torus of $SO_{2n+1}(q)$. So 
 %By Lemma \ref{hh6}, 
 $\mu(D_w)=\mu(\tilde D_{wr})$, in particular,  $\mu(D_w)$ is a $q$-character \ii  $\mu|_{\tilde D_{wr}}$ is a $q$-character. % for a maximal torus of $\tilde D_{wr}$ of $SO_{2n+1}(q)$. 
 The latter holds \ii $wr(\mu)=\mu$ by
Lemma  \ref{q-sp}. As in (1), the number of such characters equals $1_{\tilde W_j}^{\tilde W}(wr)$, but 
a priori  $h\om_j$ for $h\in \tilde W$ with $wr(h\om_j)=h\om_j$  cannot be of the shape $g\om_j$ for $g\in W$.
However, this  is the case if $j<n$, as then the orbits $W\om_j$ and $\tilde W\om_j$ coincide, so  the result follows.

 \medskip

Remark.
 The fundamental weights $\lam_j$ of the root system of type  $D_n$ are $\ep_1+\cdots +\ep_j$ for $1=1\ld j$ provided $j<n-1$, and
 $\lam_{n-1}=\frac{1}{2}(\ep_1+\cdots +\ep_{n-1}-\ep_n)$, $\lam_{n}=\frac{1}{2}(\ep_1+\cdots +\ep_{n-1}+\ep_n)$
  \cite[Planche IV]{Bo}. Therefore, $\om_n=2\lam_n$ and $\om_{n-1}=\lam_{n-1}+\lam_n$. The above results for $\GC=D_n$ can be extended to $\lam_{n-1}$. If 
  $\mu=g(\lam_{n-1}) $ $(g\in W_\GC)$ and $G=SO^+_{2n}(q)$
  then $\mu(D_w)\not\subset F_q$, whereas if $G=SO^-_{2n}(q)$ then  $\mu(D_w)\subset F_q$ \ii  
  $wg(\mu)=g(\mu)$ for $w\in W_\GC$. We do not consider this case in detail.

\begin{theo}\label{w11}  Let $G\in \{GL_n(q),SL_{n+1}(q),Sp_{2n}(q),  Spin^\pm_{2n}(q), $ q even, $SO^\pm_{2n}(q)$, q odd, $SO_{2n+1}(q), $ q odd$\}$ and let W be the Weyl group of $\GC$. Set $\om_j=\ep_1+\cdots +\ep_j$ for $1\leq j\leq n$, and assume $j<n$ if $G=SO^-_{2n}(q)$ or $Spin^-_{2n}(q)$. 
%If $G=SO^-_{2n}(q)$, assume $j\neq n$. 
Let $ \LC_j$ be a Levi subgroup of $G$ whose Weyl group is $W_j=C_W(\om_j)$, and $L_j=\LC_j^F$. Then $s_{(q-1)\om_j}$ is $L_j$-controlled. 
\end{theo}

Proof. Let $T$ be a maximal torus of $G$ and $\mu=g(\om_j)$ for $g\in W$. The weight  $(q-1)g(\om_j)$
 is trivial on $T$  \ii $g(\om_j)|_T$
is a $q$-character of $T$  (Lemma \ref{t11}).  So $(s_{(q-1)\om_j}|_{T},1_{T})$ is equal to the number of $q$-characters of $T$ in the set $g(\om_j)|_T$.
By Lemma  \ref{n1n} for $G=GL_n(q)$  and $SL_{n+1}(q)$ and by Lemma \ref{hb4}  for other classical groups above but $SO^-_{2n}(q)$, $q$ odd and $Spin^-_{2n}(q)$, $q$ even, this number is exactly $1_{W_j}^W(w)$.
So   the result follows from Lemmas  \ref{nt2} in this case. %f $G\neq SO^-_{2n}(q), Spin^-_{2n}(q)$. 
Let   
$G=SO^-_{2n}(q)$, $q$ odd, or $Spin^-_{2n}(q)$, $q$ even.
Then the number in question equals   $1_{\tilde W_j}^{\tilde W}(wr)$ in notation of Lemma \ref{hb4}(2). 
As $F_1$ arises from the graph \au which fixes the nodes $1\ld n-2$ and permutes $n-1,n$, it follows that $F_1\in C_{\tilde W}(\om_j)$ for $j=1\ld n-1$. (See the remark prior Theorem \ref{w11}.) But then, as $\tilde W=W\cdot \lan F_1\ran$, we have $\tilde W_j=C_{\tilde W}(\om_j)=C_W(\om_j)\cdot \lan F_1\ran=W_j\cdot \lan F_1\ran$. Now the result follows from Lemma \ref{nt2}, where $\tilde W_{\LC}$ is defined as $W_{\LC}\cdot \lan F_1\ran$ and $W_{\LC}$ is the Weyl group of a  Levi subgroup $\LC$.

\medskip
{\bf Proof of Theorem} \ref{th2}. 
By  Theorem \ref{w11},  the function $s_{(q-1)\om_j}$ is $L_j$-controlled.
So the result follows from Theorem \ref{dd3}.

\section{Truncated polynomials and the natural permutation module}

\subsection{The natural permutation module}

Let $G=GL_n(q)$, $V$ the natural $F_qG$-module  and  $M$  the permutation module over the complex numbers associated with the action of $G$ on non-zero vectors  of  $V$. Let $\chi$ be the character of $M$. Let $T$ be a maximal torus of $G$. 
Observe that $(\chi|_T,1_T)$   equals the number of orbits of $T$ on non-zero elements of $V$, see \cite[Theorem 32.3]{CR1962}.   We write $C_r^j$ for the $j$-th binomial coefficient (often denoted by  ${\tiny \begin{pmatrix}k\cr j\end{pmatrix}})$.

\bl{pm1} %\mar{pm1} 
Let $G=GL_n(q)$ or $SL_n(q)$ and V the natural $F_qG$-module. Let T be a maximal torus in $G$ and let  k be the composition length of T on V. 

$(1)$ Let $G=GL_n(q)$. Then the number of
T-orbits on the  non-zero vectors of $V$ is equal to $\sum_{j=1}^{k}C_k^j =
2^k-1$.

$(2)$ Let $G=SL_n(q)$. Then the  number of
T-orbits on %the  non-zero vectors of
 $V\setminus \{0\}$ is equal to %$q-2+\sum_{j=1}^{k} C_k^j =
$
q-3+2^k$.
\el

Proof. (1) Suppose first that $k=1$. Then $|T|=q^n-1$. It is well known that $T$ acts transitively on the non-zero vectors of $V$.   In general, let $V=V_1+\cdots +V_k$ be a decomposition of $V$ as a direct sum of \ir $F_qT$-modules.  Then $T=T_1\times \cdots \times T_k$, where $T_i$ acts on $V_i$ and  $|T_i|=q^{\dim V_i}-1$ ($1\leq i\leq k$).  So $T_i$ acts transitively on the non-zero vectors of $V_i$.  For $v\in V$ let  $v=v_1+\cdots+ v_k$ with $v_i\in V_i$ ($1\leq i\leq k$). 
Then $Tv=T_1v_1+\cdots +T_kv_k$. If $v_i\neq 0$  then  $T_iv_i$  consists of all non-zero vectors of $V_i$. %Therefore  
%the $T$-orbits are within $V_i$ and also if $v=v_1+\cdots+ v_k$ with $v_i\in V_i$ then the orbit of $v$ is the union of the orbits of $v_i$, except when $v_i=0$.  
%  the $T$-orbit of $v$ is the union of the $T$-orbits of $v_i$ with   $v_i\neq 0$. 
Moreover if $v'=v_1'+\cdots +v_k'$ with $v_i'\in V_i$ then  $Tv=Tv'$ \ii   $v_i=0$ implies $v_i'=0$ and conversely.  Therefore,
 the number of  $T$-orbits on $V\setminus 0$ equals   $\sum_{i=1}^k C_k^j $.

 (2) Set $T'=T\cap SL_n(q)$. Then $V_1\ld V_k$ are \ir $F_qT'$-modules so the composition length of $T'$ on $V$
  again equals $k$. Note that $|T:T'|=q-1$. Let $v=v_1+\cdots +v_k\in V$ as above. %Then $T'v=T'v_1+\cdots +T'v_k$.  
  If  $v_i=0$  for some $i$  then, obviously,   $T'v=Tv$. % coincides with $Tv$. 
     So the number of  $T'$-orbits with  $v_i=0$ for some $i$ is $\sum_{j=1}^{k-1}C_k^j $.
    
  Suppose that $v_i\neq 0 $ for all $i$. 
Let  $v'=v_1'+\cdots +v_k'$ with $0\neq v_i'\in V_i$.  Then there is a unique $t\in T$ such that $tv=v'$, so
the orbit 
 $Tv$ is regular. As $|T:T'|=q-1$, 
 the total number of  $T'$-orbits on the set $Tv$ equals $q-1$. This implies the result.

  \begin{corol}\label{c22} %\mar{c22} %Let $\chi$ be the character of $\pi$. Then
  %$(1)$ 
   $(\chi|_T,1_T)=2^k-1$ and $(\chi|_{T'},1_{T'})=q-3+2^k=q-2+(\chi|_T,1_T)$.
  \end{corol}

\bl{pm2}%\mar{pm2} 
Let $w\in S_n$, and let k be the number of cycles in the cycle decomposition of w. 
%$G=GL_n(q)$ and let $T=T_w$ be a maximal torus in $G$ and k the composition length of T on V. 
Then $2^k-1=\sum_{j=1}^{k}C_k^j =\sum_{i=1}^n 1_{Y_i}^{S_n}(w)$,
where $Y_i=S_i\times S_{n-i}$.
\el

Proof.  We know (see Proposition \ref{n1n})  that  $1_{Y_i}^{S_n}(w)$ equals the number $x_i$ of subsets $\{m_1\ld m_i\}$ of $\{1\ld n \}$
such that $w(\{m_1\ld m_i\})=\{m_1\ld m_i\}$. So  $ \sum_{i=1}^n 1_{Y_i}^{S_n}(w) $ is the sum of $x_i$,
and hence equals the number of subsets  $\{m_1\ld m_i\}$ of $\{1\ld n \}$
such that $w(\{m_1\ld m_i\})=\{m_1\ld m_i\}$, where $i$ is in the range $1\leq i\leq n$. 
% Up to conjugacy, $w$ is the product of cycles $(1\ld n_1), (n_1+1\ld n_1+n_2)\ld (n_1+\cdots +n_{k-1}\ld n)$. 
 One observes that $w(\{m_1\ld m_i\})=\{m_1\ld m_i\}$ is equivalent to saying that $\{m_1\ld m_i\}$ is a union of some cycles of $w.$
 This number can be counted in a different way, as the number of sets   $\{m_1\ld m_i\}$ for $1\leq i\leq n$ that are obtained as the union of some cycles of $w$. The number  of sets   $\{m_1\ld m_i\}$ for $1\leq i\leq n$ that consist of a single cycle is $k$,
 the number  of sets   $\{m_1\ld m_i\}$ that are unions of  two cycles is $C_k^2$,
 and  the number  of sets   $\{m_1\ld m_i\}$  that are unions of  $j$ cycles is $ C_k^j 
 $. So the result follows.

 \bl{pm34} %\mar{pm34} 
 Let $G=GL_n(q)$  and let $T$ be a maximal  torus in $G$.   
  Then $(\chi|_T,1_T)=\sum_{i=1}^n (s_{(q-1)\om_i}|_T,1_T)$, and hence
 $u(\chi\cdot St)=u(\sum_{i=1}^n (s_{(q-1)\om_i}\cdot St)$. 
 \el

Proof.  Let  $T=T_w$ for $w\in W\cong S_n$ and let $k$ be the composition length of $T$ on $V$.
% (or the number of cycles in the cycle decomposition of $w$). 
Then   $(\chi|_T,1_T)=2^k-1 %\sum_{j=1}^{k} C_k^j 
 $ by Lemma \ref{pm1}, and $(s_{(q-1)\om_i}|_T,1_T)=1_{Y_i}^{S_n}(w)$ by Theorem \ref{w11}.
 So the result follows from Lemma \ref{pm2}. The second  statement follows from Lemma \ref{ed2}.

\bl{pm4} %\mar{pm4} 
Let $G=GL_n(q)$ and let $\nu$ be a unipotent character  of $G$.
 % Let $\chi$ be the character of M. 
  Then   $(\nu,\chi\cdot St)=(\nu,\sum_{i=1}^n s_{(q-1)\om_i}\cdot St)=(\nu,\sum_{i=1}^n St_i^{\# G})$, where 
 $St_i$ is the Steinberg character of a Levi subgroup of $G$ isomorphic to $GL(i)\times GL(n-i,q)$.  
 \el

Proof. The former equality  is Lemma  %s \ref{dd3} and
 \ref{pm34},  the latter one follows from Theorem \ref{th2}. %s \ref{w11} and \ref{dd3}. 
 Note that the second equality can also be deduced   from \cite[Theorem 6.2]{HZ1}.

\bl{pp3} %\mar{pp3} 
Let $G=GL_n(q)$, $G'=SL_n(q)$ and let V be the natural $F_qG$-module. Let %$\chi$ be as above and 
 $\phi=\chi|_{G'}$. Let $St$, $St'$ be the Steinberg characters of $G,G'$, respectively. Then   $(St,\chi\cdot St)=n$ and   $(St',\phi\cdot St')=n+q-2$.
\el

Proof.  It is well known that $(St, St_i^{\# G})=1$. (Indeed, if $\si$ denotes the Harish-Chandra restriction of $St$ to $L_i$ then $\si=St_i$ by  \cite[p.72]{DM}. By Harish-Chandra reciprocity  \cite[Proposition 70.(iii)]{CR-2},  $(St, St_i^{\# G})=(\si, St_i)=1$.)   So Lemma \ref{pm4} yields the result for $G$.

The Weyl groups of  $\GC$ and $\GC'$ coincide. So the mapping  $\TC \ra \TC'=\TC\cap  \GC'$
yields a bijection between the $G$-conjugacy classes of $F$-stable maximal tori in $\GC$ 
and the  $G'$-conjugacy classes of $F$-stable maximal tori in $\GC'$.
In addition, $W(T)=W(T')$ and     $(\phi|_{T'},1_{T'})=q-2+
(\chi|_T,1_T)$ by  Corollary \ref{c22}. So formula (\ref{st4}) yields

$$(St', \phi\cdot St')= \sum _{\TC'  }
\frac{ (\phi|_{T'},1_{T'})}{|W(T')|}= \sum _{\TC' }
\frac{ q-2+(\chi|_T,1_T)}{|W(T')|}
= \sum _{\TC  }
\frac{ (\chi|_{T},1_{T})}{|W(T)|}+(q-2)\sum _{\TC' }
\frac{ 1}{|W(T')|}.$$
By formula (\ref{st4}) applied to $GL_n(q)$, we have $ \sum _{\TC  }
\frac{ (\chi|_{T},1_{T})}{|W(T)|}= (St,\chi\cdot St)$.
The latter equals $n$. By formula (\ref{eq2}),  $1=(St',St')=\sum _{\TC'  }
\frac{ 1}{|W(T')|}$, whence  the result.

\subsection{Truncated polynomials}

Let $R_n=\overline{F}_q[X_1\ld X_n]$ be the \po ring with indeterminates $X_1\ld X_n$ over $\overline{F}_q$.
Let $I$ be the ideal of $R_n$ generated by $X_1^p\ld  X_n^p$ and $\overline{R}_n=R_n/I$.
Then $\overline{R}_n$ can be viewed as the truncated \po ring whose elements are linear combinations of monomials $X_1^{c_1}\cdots X_n^{c_n}$ with  $0\leq c_1\ld c_n<p$. Let $\GC=GL_{n}(\overline{F}_q)$, $G=GL_n(q)$. Viewing $X_1\ld X_n$ as a basis of the natural $\overline{F}_q\GC$-module ${\mathbf V}_n$, one extends the action of $\GC$ to $\overline{R}_n$ in a standard way. Note that homogeneous \pos of a fixed degree  form an $\overline{F}_q\GC$-submodule 
of  $\overline{R}_n$. 
  We often view  $\overline{R}_n$ as an $\overline{F}_qG$-module
and as an $\overline{F}_q\GC'$-module, where $\GC'=SL_{n}(\overline{F}_q)$.

Let $q=p^m$. There is a well known embedding $G=GL_n(q)\ra GL_{mn}(p)$ obtained by viewing $F_q$ as a vector space over $F_p$. So the above constructed $F_pGL_{mn}(p)$-module $\overline{R}_{mn}$ 
is also an $F_pGL_n(q)$-module, which plays a significant role below.

The action of $G$ on $\overline{R}_{mn}$  extends to $\GC$ as follows. For $g\in \GC$
consider the mapping $g\ra \diag(g,Fr(g)\ld Fr^{m-1}(g))$,
 $Fr(g)$ is obtained from $g$ by raising every matrix entry to the $p$-power and  $\diag (g,Fr(g)\ld Fr^{m-1}(g))$ means the  block diagonal matrix with diagonal blocks $g,Fr(g)\ldots $. This makes the space 
${\mathbf V}_{mn}$ to be a $\GC$-module. One easily observes that the two actions of $G$ on 
${\mathbf V}_{mn}$ obtained from the embeddings $G\ra GL_n(\overline{F}_{q})$ and $G\ra GL_{mn}(p)$
yield equivalent \reps of $G$. 
Note that $\overline{R}_{mn}$ is completely reducible both as $\overline{F}_{q}\GC$- and $\overline{F}_{q}G$-module (see \cite[Proposition 1.6]{SZ2}).  

Denote by $M_q$ the $\overline{F}_{q}G$-module obtained from the action of $G$ on the vectors 
of the natural $F_qG$-module $V$ (the zero vector is not excluded). The Brauer character of $M_q$
coincides with $\chi+1_G$ on the $p$-regular elements. 

A remarkable fact going back to Bhattacharia \cite{Bh} states that the Brauer characters of  $G$ on $M_q$ and on  
$\overline{R}_{mn}$ are the same. This was exploited in  \cite{SZ3} to obtain the decomposition numbers of the \ir consituents of $M$. The following lemma is essentially  \cite[Theorem 3.2]{SZ3}.

\bl{bh6} %\mar{bh6} 
Let $G=GL_n(q)$. The Brauer characters  of $\overline{R}_{mn}|_G$ and $M_q$ coincide. 
\el

Proof. Let $V$ be the natural $F_qG$-module and  $V_m$  the natural $F_pGL_{mn}(p)$-module. Viewing $V$ as a vector space over $F_p$, we identify the additive group of $V$ with that of $V_m$. Moreover, the regular embedding $F_q\ra Mat(m,F_p)$ yields an embedding  $h:GL_n(q)\ra GL_{mn}(p)$, and the permutation actions 
of $G$ on $V$ and $V_{m}$ are isomorphic. Therefore, the permutation characters afforded by these actions coincide. For $m=1$ (that is for $q=p$) the lemma is proved by Bhattacharia \cite{Bh}. Clearly, this remains true for every subgroup of $GL_n(p)$. Applying this to the subgroup $h(G)$ of $GL_{mn}(p)$, we obtain the result for $G$.   

\medskip
 Denote  $\lam_1\ld \lam_{n-1}$  the fundamental weights of $\GC'=SL_n(\overline{F}_q)$. 
 Weights $a_1\lam_1+\cdots +a_{n-1}\lam_{n-1}$
with $0\leq a_1\ld a_{n-1}< q $  are called $q$-restricted, and 
 the mapping $\rho_\lam\ra \rho_\lam|_{G'}$ sets up a bijection between the \ir \reps of $\GC$ with $q$-restricted highest weights  and \ir \reps of $G'=SL_n(q)$ over 
 $\overline{F}_{q}G$.

\bl{d1d} % \mar{d1d} 
Let $\lam=a_1\lam_1+\cdots +a_{n-1}\lam_{n-1}$ be a weight of $\GC'=SL_n(\overline{F}_{q})$. 

$(1)$ If $\lam$ occurs as a weight of the $\GC'$-module $\overline{R}_{n}$  then $-p<a_i<p$ for $i=1\ld n-1$. 

$(2)$ Suppose that  $\lam$ is dominant. Then  $\lam$ is a weight of the $\GC'$-module $\overline{R}_{n}$ \ii $\lam$ is strongly $p$-restricted, that is, $a_1+\cdots +a_{n-1}< p$.
\el

Proof. Let $\TC'$ be the group of diagonal matrices in $\GC'$. Then monomial polynomials are weight vectors for $\TC'$ and they form a basis of $\overline{R}_{n}$. Let $f=X_1^{c_1}\cdots X_n^{c_n}\in \overline{R}_{n}$, so $0\leq c_1\ld c_{n}< p$.
Then the weight of $f$ for $\GC'$ in terms of $\ep_1\ld \ep_n$ is $c_1\ep_1+\cdots +c_n \ep_n$. Note that $\ep_1+\cdots + \ep_n$ is the zero weight for $\GC'$, so we can write $\ep_n=-(\ep_1+\cdots+ \ep_{n-1})=-\lam_{n-1}$.
 As 
$\ep_1=\lam_1$  and $\ep_i=\lam_i-\lam_{i-1 }$
for $1<i<n$, the weight of $f$ in terms of 
$\lam_1\ld \lam_{n-1}$ is $(c_1-c_2)\lam_1+\cdots +(c_{n-1}-c_n)\lam_{n-1}$. Set $a_i=c_i-c_{i+1}$ for $i=1\ld n-1$. As $0\leq c_i<p$, (1) follows. 

(2)  Suppose that  $\lam$ is a dominant  weight of  $\overline{R}_{n}$.  Then  $a_i= c_i-c_{i+1}\geq0$.    This is equivalent to   the condition  $c_1\geq \cdots \geq c_n\geq 0$.  Furthermore,  $a_1+\cdots +a_{n-1}=c_1-c_n< p$.
% so $a_1+\cdots +a_{n-1}\leq p-1$. 
  Conversely, suppose that $a_1+\cdots +a_{n-1}< p$. Set $c_i=a_i+\cdots +a_{n-1}$ for $i<n$ and $c_n=0$. Then the monomial $X_1^{c_1}\cdots X_{n-1}^{c_{n-1}}$ has weight $\lam$ and $c_i< p$ for $i=1\ld n$. 

\medskip 
Remark. Let $\al_0$ be the longest root of $\GC'$. Then $a_1+\cdots +a_{n-1}=\lan \lam,\al^{\rm v}_0\ran$
in notation of  \cite[p.16]{Hu}.    

\bl{zz1} %\mar{zz1}
 $(1)$  Let $\GC=GL_n(\overline{F}_{q})$.  There is a $ \GC$-module isomorphism $\overline{R}_{mn}\cong \overline{R}_{n}\otimes Fr (\overline{R}_n)\otimes\cdots\otimes Fr^{m-1}(\overline{R}_n) $.  
 
 $(2)$ Let $\GC'=SL_n(\overline{F}_{q})$.  Then every weight  $\mu$ of $\GC'$ on $\overline{R}_{mn}$ can be expressed as  $\mu_0+\mu_1p+\cdots +\mu_{m-1}p^{m-1}$, where $\mu_0,\ld \mu_{m-1}$ are weights of $\GC'$ on  $\overline{R}_{n}$.
 
 $(3)$ If $\mu=0$ or $(q-1)\lam_k$  $(1\leq k\leq n-1)$ then the expression in $(2)$ is unique. 
\el

Proof. (1) See \cite[Proposition 1.6]{SZ2}. Observe that $Fr (\overline{R}_n)$ is the $\GC$-module
obtaining from $\overline{R}_n$ by twisting with morphism $Fr:\GC\ra \GC$ defined in the beginning 
of the section (so $g\in \GC$ acts via $Fr(g)$).  (2)  follows from (1).  

(3) Let $\mu_i=a_{i,1}\lam_1+\cdots +a_{i,n-1}\lam_{n-1}$ $(0\leq i<m)$
so $\mu=\sum_{i=0}^{m-1}\sum_{j=1}^{n-1} a_{ij}p^i\lam_j=\sum_{j=1}^{n-1}(\sum_{i=0}^{m-1} a_{ij}p^i)\lam_j$. 
By Lemma \ref{d1d}, $-p<a_{i,j}<p$ for $j=1\ld n-1$. Let $\mu=0$.  Then $\sum_{i=0}^{m-1} a_{ij}p^i=0$ for every $j$.
If  $\mu_i\neq 0$ for some $i$, then $a_{ij}\neq 0$ for some $j$.  %Fix this $j$.
Let $r=\max\{ i: a_{ij}\neq 0\}$. We can assume $a_{rj}>0$. Then
 $0=\sum_{i=0}^{m-1} a_{ij}p^i= a_{rj}p^r+\sum_{i=0}^{r-1} a_{ij}p^i\geq p^r-(p-1)(1+p+\cdots +p^{r-1})=1$, which is a contradiction. 
 
 Let $\mu= (q-1)\lam_k$. 
 Suppose that $\mu=\mu_0+\mu_1p+\cdots +\mu_{m-1}p^{m-1}$, where $\mu_0,\ld \mu_{m-1}$ are weights of  $\overline{R}_{n}$. Then  $\sum_{j=1}^{n-1}(\sum_{i=0}^{m-1} a_{ij}p^i)\lam_j=(q-1)\lam_k$, whence $\sum_{i=0}^{m-1} a_{ij}p^i=0$ for $j\neq k$. We have seen in the previous paragraph that this implies $a_{ij}=0$ for $i\neq k$. 
%So $\mu_j=a_{kj}\lam_k$. 
Then  $(q-1)\lam_k=\mu=(a_{k0}+a_{k1}p+\cdots +a_{k,m-1}p^{m-1})\lam_k$, and hence   $q-1=a_{k0}+a_{k1}p+\cdots +a_{k,m-1}p^{m-1}\geq (p-1)(1+p+\cdots +p^{m-1})=q-1$. The equality holds \ii $a_{k0}=a_{k1}=\cdots =a_{k,m-1}=p-1$. So the result follows.

 \bl{bbb} %\mar{bbb} 
 Let $\lam$ be a strongly $q$-restricted  weight of $\GC'$ and $\rho_\lam$ an \irr of $\GC'$ with highest weight $\lam$. 
 Then all weights of $\rho_\lam$ occur in $\overline{R}_{mn}.$
 \el

 Proof. Suppose first that $m=1$. By Lemma \ref{d1d}, $\lam$ is a weight of $\overline{R}_{n}.$ Let $\tau$ be an  \irr of $\GC'$ on  $\overline{R}_{n}$ such that $\lam$ is a weight of $\tau$. Then the weights of $\rho_\lam$ are weights of $\tau$ (but $\rho_\lam$ is not usually a constituent of  $\overline{R}_{n}$.) This follows from a criterion in \cite[Ch.VIII, Corollary 2 of Proposition 3]{Bo}  for \reps in characteristic 0, and from Suprunenko's theorem  \cite{Su1} or
 \cite[\S 3.3]{Hu} for prime characteristics.

 Let $m>1$.    
 Then $\lam=\nu_0+p\nu_1+\cdots+p^{m-1}\nu_{m-1}$, where  $\nu _0\ld \nu_{m-1}$ are strongly $p$-restricted weights. By Lemma \ref{d1d}, $\nu_0\ld \nu_{m-1}$ are weights of $\overline{R}_{n}$. By (1),  all weights of $\rho_{\nu_i}$
 $(i=0\ld m-1)$ occur  in $\overline{R}_{n}$. As $\rho_\lam=\rho_{\nu_0}\otimes Fr(\rho_{\nu_1})\otimes \cdots \otimes Fr^{m-1}(\rho_{\nu_{m-1}})$,   all weights of $\rho_\lam$ occur in $\overline{R}_{mn}$.

\medskip
Remark. If $q>p$ then not every weight of $\overline{R}_{mn}$ is strongly $q$-restricted. 
For instance let $p=n=3$, $q=9$ and $\lam=2\lam_1+2\lam_2$. As $\lam=-(\lam_1+\lam_2)+3(\lam_1+\lam_2)$
and $-(\lam_1+\lam_2)$ is a weight of $\rho_{\lam_1+\lam_2}$, it follows that $\lam$ is  a  weight of 
the \rep $\rho_{\lam_1+\lam_2}\otimes Fr(\rho_{\lam_1+\lam_2})$, which is  a constituent of $\overline{R}_{mn}$
for $m=2,n=3$.

\bl{zw1} %\mar{zw1}
 $(1)$ Let $1\leq k\leq p-1$. The polynomial $f_k=(X_1\cdots X_n)^{k}\in \overline{R}_n$ is a vector of zero weight for  $\GC'$. 

 $(2)$  Let  $f=X_1^{a_1}\cdots X_n^{a_n}\in \overline{R}_n$,  let $e=\max_{1\leq i,j\leq n}(a_i-a_j)$ and let $\nu$ be the $\GC'$-weight of $f$. Then the \mult of  $\nu$ in  $\overline{R}_n$ is equal to $p-e$.  
 
 $(3)$  The \mult of weight $(q-1)\lam_i$ in $ \overline{R}_{mn}$ equals $1$ and that of weight $0$
equals $q$. 
\el

Proof.  Let $\TC'$ be the group of diagonal matrices of $SL_n(\overline{F}_{q})$. Note that monomials $X_1^{a_1}\cdots X_n^{a_n}$ form a basis of $ \overline{R}_n$, and each monomial is a weight vector.  

 $(1)$   Let $t=\diag(t_1\ld t_n)\in T'\subset SL(n,\overline{F}_{q})$. Then $f_k(t)=(\det t)^k=1$.

  $(2)$      If all $a_1\ld a_n$ are non-zero then there is $k>0$ such that $ f=f_k\cdot  X_1^{a_1-k}\cdots X_n^{a_n-k}$ and some $a_i-k=0$. 
Clearly, the weight of $f':=X_1^{a_1-k}\cdots X_n^{a_n-k}$ is the same as that of $f$, so we can assume 
that $f=f'$, that is, that $a_j=0$ for some $j\in \{1\ld n\}$. Let $h= X_1^{b_1}\cdots X_n^{b_n}$ be another monomial of weight $\nu$ and  $b_l=0$ for some $l\in\{1\ld n\} $. 
Observe that $a_j=0$ \ii $b_j=0$. Indeed, suppose the contrary, say, let $a_j=0$, $b_j\neq0$,  and $b_l=0$ for some $l>j$. Let $t=\diag (1\ld 1,t_j,1\ld 1,t_j\up,1\ld 1)\in SL(n,\overline{F}_{q})$, where $0\neq t_j\in \overline{F}_{q}$
is arbitrary and $t_j\up$ occupies the $l$-th position. Then
we have $t\cdot f= t_j^{-a_l}\cdot f$ and $t\cdot h=t_j^{b_j}\cdot h $. As $f,h$ are weight vectors of the same weight $\nu$, we have  $ t_j^{-a_l}= t_j^{b_j}$. 
As $a_l\geq 0, b_j\neq 0$, this is a contradiction. Similarly, $b_j=0$ implies $a_j=0$. Finally, fixing $j$ with 
$a_j=b_j=0$, let $t=\diag(t_1\ld t_n)$ with $t_j=(\Pi_{i\neq j}t_i)\up $. One easily observes that $t\cdot f=
(\Pi _{i\neq j}t_i^{a_i})f=\nu(t)f$ and  $t\cdot h=
(\Pi _{i\neq j}t_i^{b_i})h=\nu(t)h$. As each $t_i$ with $i\neq j$ is arbitrary, $a_i=b_i$ for all $i\neq j$.

\itf if $f=X_1^{a_1}\cdots X_n^{a_n}$ is of weight $\nu$ with $a_j=0$ for  some $j$ then all other monomials of weight $\nu$ must be 
$f\cdot f_k$ for some $k$.  In this case
 $e=\max_ia_i$.  Therefore, $f\cdot f_k\in \overline{R}_n$ implies $e+k\leq p-1$, so the result follows.

 $(3)$  The monomial $X_1^{p-1}\cdots X_i^{p-1}\in \overline{R}_{n}$ is a vector of  weight $(p-1)\lam_i$
 and $f_1$ is that of weight 0. So for $m=1$ the lemma follows by applying (2) to these  monomials. 
Since  $
\overline{R}_{mn}\cong \overline{R}_{n}\otimes Fr (\overline{R}_{n})\otimes\cdots \otimes Fr^{m-1}(\overline{R}_{n})$,  
the result for arbitrary $m$ follows from Lemma 
 Lemma \ref{zz1}(3). 

\bl{m1m} %\mar{m1m}
 Let  $\nu$ be a $\GC'$-weight on $\overline{R}_{mn}$ and T  a maximal torus of $G'=SL_n(q)$. Suppose that $\nu|_T=1_T$.
Then either $\nu=0$ or $\nu$ lies in the $W$-orbit of $(q-1)\lam_i$ for some $i\in \{1\leq i\leq n-1\}$.
\el

Proof.   Note that
 $\overline{R}_{mn}$ is the sum of weight spaces of  $\TC$ on $\overline{R}_{mn}$, so
$(\overline{R}_{mn}|_T,1_T)=\sum _\mu m_\nu (s_\mu|_T,1_T)$, where $\mu$ runs over the dominant weights of  $\GC'$ in $\overline{R}_{mn}$ and $m_\nu$ is the \mult of $\nu$. By Lemma \ref{zw1}, the \mult of  weight 0  of $\GC'$ in  $\overline{R}_{mn}$ equals $q$ and that of $(q-1)\lam_i$ is 1 for $i\in \{1\ld n-1\}$. Therefore, $(\overline{R}_{mn}|_T,1_T)=q+\sum_{i=1}^{n-1}(s_{(q-1)\lam_i}|_T,1_T)+\sum _{\mu } m_\nu (s_\mu|_T,1_T)$, where $\mu$ now runs over the dominant weights of  $\GC'$ in $\overline{R}_{mn}$ distinct from 0 and $(q-1)\lam_i$ for $1\leq i<n $. Let $T=T_w$ for $w\in W$. By Proposition \ref{n1n} (together with Lemma \ref{t11}),  $(s_{(q-1)\lam_i}|_{T_w}, 1_{T_w})= 1^W_{W_i}(w)$. 
By Lemma \ref{pm2},
$\sum_{i=1}^{n-1}1^W_{W_i}(w)=2^k-2$, where $k$ is the number of cycles in the cycle decomposition of $w$ (as the term with $i=n$ in Lemma \ref{pm2} equals 1). So 
$(\overline{R}_{mn}|_T,1_T)=q+2^k-2
+\sum 
m_\nu (s_\mu|_T,1_T)$, where $\mu\neq 0, (q-1)\lam_1\ld (q-1)\lam_{n-1}$.

Another calculation (Corollary  \ref{c22} and Lemma \ref{bh6}) shows that $(\overline{R}_{mn}|_T,1_T)=q-2+2^k$. \itf  $\sum_\mu (s_\mu|_T,1_T)=0$ for every dominant weight $\mu$ of $\overline{R}_{mn}$
distinct from 0 and $(q-1)\lam_i$ for $i\in \{1\ld n\}$. This implies the lemma.  

 \medskip
 Remark. At the first sight, Lemma \ref{m1m} cannot be true for $q=2$ since the split torus $T$ has order 1, and hence $\nu|_T=1_T$ for arbitrary weight $\nu$. However, in fact the Brauer character of $\overline{R}_{n}$
 equals $2\cdot 1_T+\sum_{i=1}^{n-1}s_{\lam_i}$, which agrees with the statement of the lemma.

\bl{bh2} %\mar{bh2} 
Let $j\in\{1\ld n-1\}$,  $\GC=GL_n(\overline{F}_{q})$ and ${\mathbf G}'=SL(n,\overline{F}_{q})$. 
 Then ${\mathbf G}'$-module $\overline{R}_{mn}$ contains composition factors of highest weights $(q-1)\lam_j$
$(j=1\ld n-1 )$, 
each  occurs with \mult $1$.  If $\GC=GL_n(\overline{F}_{q})$  
then the $\GC$-module $\overline{R}_{mn}$ contains a composition  factor of highest weight $(q-1)\om_j$   with \mult $1$.
 \el

Proof. Suppose first that $m=1$. Then $v_j:=X_1^{p-1}\cdots X_j^{p-1}$ is a vector of weight $(p-1)\lam_j$ for $\GC'$ and of weight $(p-1)\om_j$ for $\GC$. Moreover, $gv_j=v_j$ for every upper unitriangular matrix
in $GL_n(\overline{F}_{q})$.  As the space of homogeneous \pos of fixed degree is an \ir $\GC'$-module \cite{Bh,SZ2},   $v_j$ is a highest weight vector in this module (both for $\GC'$ and $\GC$). By Lemma \ref{zw1}, the \rep afforded by this module occurs in  $\overline{R}_{n}$ with \mult 1. 

As  $\overline{R}_{mn}\cong  \overline{R}_{n}\otimes Fr( \overline{R}_{n})\otimes Fr^{m-1}( \overline{R}_{mn})$, it follows that  $\overline{R}_{mn}$ contains a composition factor in question, 
and the statement on \mult follows again from Lemma \ref{zw1}(3).

\medskip
If $n>2$ then Lemma \ref{bh2} is a special case of \cite[Theorem 1.4]{SZ2}, which determines the \mult of every 
 composition factors of $\GC'$ on $\overline{R}_{mn}$ for $n>2$.

   \begin{lemma}\label{cs4} %\mar{cs4} 
 Let $G'=SL_n(q)$ and $T=T_w$ a maximal torus of $G'$.   %Let $\nu\neq 0$ be a strongly q-restricted weight of $\GC'=SL_{n}(\overline{F}_q)$. %occurring in $\overline{R}_{mn}$.  
    Let $\rho_\nu$ be an \irr of $\GC'=SL_{n}(\overline{F}_q)$ with a strongly q-restricted highest weight $\nu\neq 0$, and $d_0$  the \mult of weight $0$ in $\rho_\nu$. Let $\beta_\nu$ be the Brauer character of $\rho_\nu$.
 
 $(1)$  $(\beta_\nu|_T,1_T)=d_0$  
 unless $\nu=(q-1)\lam_i$ for $i\in \{1\ld n-1\}$ and $w\neq [n]$.    
 
 $(2)$ Let $\nu=(q-1)\lam_i$ for some $i\in \{1\ld n-1\}$. Then   $(\beta_\nu|_T,1_T)$ equals  $d_0+1_{W_i}^W(w)$ and    $d_0\leq 1$,
 and  $d_0=1$ \ii  $i(p-1)\equiv 0\pmod n$.
\end{lemma}

Proof. (1) By Lemma \ref{bbb}, every weight of $\rho_\nu$ occurs as a weight of  $\overline{R}_{mn}$.
Note that $\beta_\nu|_T=d_0\cdot 1_T+\sum_\mu d_\mu s_\mu|_T$, where the sum is over 
 the dominant weights $\mu\neq0$  of $\rho_\nu$ and $d_\mu$ denotes the \mult of $\mu$  in $\rho_\nu$.   By Lemma \ref{m1m}, 
%By Lemma \ref{ze1}, 
either $(s_\mu|_T,1_T)=0$ or $\mu=(q-1)\lam_i$ for some $i\in\{1\ld n-1\}$. Suppose that the latter  holds. Then $\mu=\nu$. Indeed,  $(q-1)\lam_j$ occurs as a weight  $\overline{R}_{mn}$ with \mult 1 (Lemma \ref{m1m}(3)), and $\rho_{(q-1)\lam_j}$ is a constituent of $\overline{R}_{mn}$ (Lemma \ref{bh2}).  Therefore, $(q-1)\lam_i$ is a weight of $\rho_\nu$ \ii
$\nu=(q-1)\lam_i$. Thus, if $\nu\neq (q-1)\lam_i$ for some $i$ then $(\beta_\nu|_T,1_T)=d_0$.

 Let $\nu=(q-1)\lam_i$ for some $i\in\{1\ld n-1\}$. 
As above, $\mu\neq (q-1)\lam_j$ for $j\neq i$. 
By Lemma \ref{n1n}, $(s_{(q-1)\lam_i}|_T,1_T)=1_{W_i}^W(w)$. 
As $1_{W_i}^W(w)=0$ for  $w=[n]$, (1) follows.  

(2) We show that $d_0\leq 1$ for $\nu=(q-1)\lam_i$. Indeed,  $\rho_\nu$ is a constituent of $\overline{R}_{mn}$.
Every weight of every \ir constituent of $\overline{R}_{mn}$ has \mult 1.
This easily follows from the fact that $\chi|_{T_n}$ is the character of the regular \rep of $T_n$,
where $T_n$ is the maximal torus of $GL_n(q)$ of order $q^n-1$ (see also the proof of Lemma \ref{pm1}).
Indeed, this implies, in view of Lemma \ref{bh6}, that all weights of $\GC=GL_n(\overline{F}_{q})$ on the space of non-constant \pos of $\overline{R}_{mn}$ have \mult 1. As $\GC={\bf Z}\cdot \GC'
$, where ${\bf Z}$ is the center of   $\GC$, the weight \mult of an \irr of $\GC $ are the same as those of the restriction to $\GC'$. 

Furthermore,  Lemma \ref{zz1} tells us that $\rho_{(q-1)\lam_i}$ has  weight 0 \ii so is $\rho_{(p-1)\lam_i}$.  This happens \ii  $(p-1)\lam_i$ lies in the root lattice. 
(For Lie algebras over the complex numbers this is well known \cite[Ch. VIII, \S 7, Proposition 5(iv)]{Bo2}. By a result of Suprunenko \cite{Su1}, this remains true for simple algebraic groups of type $A$ and \ir \reps with  $p$-restricted highest  weights. See also \cite[\S 3.3]{Hu}.)  The quotient  of the weight lattice over the  root  lattice is isomorphic to $\ZZ \pmod n$, the cyclic group of order $n$,  and $i\lam_1-\lam_i$ for $i<n$ is a linear combination of roots \cite[Planche I]{Bo}. 
It follows that  $(p-1)\lam_i$ lies in the root lattice \ii $i(p-1)\equiv 0\pmod n$.

\medskip
Remark. In general,  $\rho_\nu$ needs not be a composition factor of $\overline{R}_{mn}$ if $\nu$ is a dominant weight of
$\overline{R}_{mn}$. Moreover, if $\rho_\nu$ has a weight 0, then the \mult of it equals 1 \ii $\rho_\nu$ is a
consituent of $\overline{R}_{mn}$.

\begin{corol}\label{c56} Let $G'$, $\nu$, $\beta_\nu$, $d_0$ be as in  Lemma $\ref{cs4}$. 

$(1)$ Let $w=[n]$ and $T=T_w$, so $|T|=(q^n-1)/(q-1).$ Then $(\beta_\nu|_T,1_T)=d_0$. 

$(2)$  Let $w=[n-1,1]$ and $T=T_w$, so $|T|=q^{n-1}-1.$ Then $(\beta_\nu|_T,1_T)$ equals the \mult of weight $0$ in $\rho_\nu$, unless $\nu=(q-1)\lam_1$ or $(q-1)\lam_{n-1}$. In these cases $(\beta_\nu|_T,1_T)=1+d_0$, where $d_0= 1$ if n divides $p-1$, otherwise $d_0=0$. 
\end{corol}

Proof. (1) In this case $1_{W_i}^W(w)=0$ as $W_i$ contains no conjugate of $w$ for every $i\in \{1\ld n-1\}$. So we are done by   Lemma \ref{cs4}. %Theorem \ref{th5}.

(2) In this case $W_i $ contains no conjugate of $w$ \ii $i=1$ or $n-1$. Then  $1_{W_1}^W$ is the permutation \rep of $W\cong S_n$ associated with the natural action of $S_n$ on $n$ letters. So $1_{W_1}^W(w)=1$. As $W_{n-1}=W_{1}$, the result follows from  Lemma \ref{cs4}.  

\medskip
Remark. The special case of Corollary \ref{c56}(2) for $n=2$ and $q$ even is examined in \cite[Corollary of Theorem 3]{KS}. Note that in this case
every $q$-restricted weight $\nu$ of $\GC$ is strongly $q$-restricted, and $\rho_\nu$ has weight zero only if $\nu=0$. 

We also mention  results of  \cite[Theorem 8]{zv} and \cite[Theorem 1]{KS} where the authors prove existence  of eigenvalue 1 of elements $g$ 
of cyclic tori $T_w$ with $w=[n]$ or $[n-1,1]$ in representations $\rho$   under certain restrictions on $|g|$ and  $\rho$.

 \subsection{Unipotent constituents of certain projective modules}

 It is well known that the tensor product of any $\overline{F}_qG$-module with a projective module is again a projective module. 
Let $St_q$ denote the Steinberg 
$\overline{F}_qG$-module. Then $St_q$ is a projective $\overline{F}_qG$-module. If $\rho$ is an \irr of $G'=SL_n(q)$ with Brauer character $\beta$ then $\rho\otimes St_q$ is a projective  $\overline{F}_qG$-module, whose  character is $\beta\cdot St$.  
(If $\beta$ is a generalized Brauer character  of $G$ then $\beta\cdot St$ is meant to be the function equal to 0
at the $p$-singular elements of $G$ and $\beta(g)\cdot St(g)$ for $g\in G$ semisimple.)

\begin{theo}\label{ut1} %\mar{ut1}
 Let $\GC'=SL_{n}(\overline{F}_q)$ and $G'=SL_n(q)$. Let $\rho_\nu$  be an \irr of $\GC'$ with highest weight $\nu$,
 $d_0$ the \mult of weight $0$ of $\rho_\nu$ and $\beta_\nu$ the Brauer character of $\rho_\nu|_{G'}$.    
Suppose that $\nu$ is 
strongly q-restricted.

$(1)$  $u(\beta_\nu\cdot St)=d_0\cdot St$, unless $\nu=(q-1)\lam_i$
for $i\in\{1\ld n-1\}$.

$(2)$ Suppose that  $\nu=(q-1)\lam_i$ for  $i\in\{1\ld n-1\}$, and let $L_i$ be a Levi subgroup of $G$
whose Weyl group is $W_i:=C_{W}(\lam_i)$. Then $u(\beta_\nu\cdot St)=d_0\cdot St+St_{L_i}^{\# G}$,
where $d_0\leq 1$. 

\end{theo}

Proof. (1) 
By Lemma \ref{ed2}, $\beta_\nu\cdot St=\sum_{\TC}\frac{(\beta_\nu|_T,1_T)}{|W(T)|}\ep_{\GC'}\ep_{\TC}R_{\TC,1}$, where the sum is over representatives of $G$-conjugacy classes of $F$-stable maximal tori $\TC$ of $\GC'$, and $T=\TC^F$. By Lemma \ref{cs4}, if   $\nu\neq
(q-1)\lam_i$ for some $i\in\{1\ld n-1\}$ then $(\beta_\nu|_T,1)=d_0$. Therefore, 
$\beta_\nu\cdot St=d_0\cdot \sum_{\TC}\frac{1}{|W(T)|}\ep_{\GC'}\ep_{\TC}R_{\TC,1}=d_0\cdot St$
by formula (\ref{eq2}).

(2) Let  $\nu=(q-1)\lam_i$ and $T=T_w$. Then $(\beta_\nu|_T,1)=d_0+1_{W_i}^W(w)$ and $d_0\leq1$  by Lemma \ref{cs4}(2). So
$\beta_\nu\cdot St=d_0\cdot \sum_{\TC}\frac{1}{|W(T)|}\ep_{\GC'}\ep_{\TC}R_{\TC,1}+
\sum_{\TC}\frac{1_{W_i}^W(w)}{|W(T)|}\ep_{\GC'}\ep_{\TC}R_{\TC,1}$. It follows from 
Lemma \ref{nt2}      that the function $\beta_\nu-d_0$ is $L_i$-controlled. 
The former sum equals $d_0\cdot St$ (as above), whereas the latter sum yields $St_{L_i}^{\# G'}$ 
by Theorem \ref{dd3}. 

 \medskip
 Remark. The condition for $d_0=1$ is given in Lemma \ref{cs4}(2).

\medskip
{\bf  Proof of Theorem} \ref{th5}.
The results stated follow from Lemma \ref{cs4} and Theorem \ref{ut1}.

\end{document}